\documentclass[12pt]{article}
\usepackage{amssymb,amsmath,amsthm}

\title{\bf \normalsize Nguetseng's  Two-scale Convergence  Method\\
 For Filtration and Seismic Acoustic Problems\\
in Elastic Porous  Media}

\author{Anvarbek~Meirmanov\footnote{This work is partially supported by the grant of Higher
Education Commission of Pakistan under National Research Program
for Universities. The title of the grant is:
\textit{Homogenization of Underground Flows} }}
\date{}

\theoremstyle{plain}
\newtheorem{theorem}{Theorem}[section]
\newtheorem{lemma}{Lemma}[section]

\newtheorem{corollary}{Corollary}[section]

\theoremstyle{definition}
\newtheorem{definition}{Definition}[section]
\newtheorem{assumption}{Assumption}[section]

\theoremstyle{remark}
\newtheorem{remark}[theorem]{Remark}

\numberwithin{equation}{section}
\renewcommand{\div}{\mbox{div}}
\newcommand{\x}{{\mathbf x}}
\newcommand{\y}{{\mathbf y}}
\newcommand{\w}{{\mathbf w}}
\newcommand{\h}{{\mathbf h}}
\newcommand{\uu}{{\mathbf u}}
\newcommand{\vv}{{\mathbf v}}
\newcommand{\e}{{\mathbf e}}
\newcommand{\n}{{\mathbf n}}
\renewcommand{\mathbf}[1]{\mbox{\boldmath$#1$}}
\newcommand{\V}{\mathbf V}
\newcommand{\W}{\mathbf W}
\newcommand{\UU}{\mathbf U}
\newcommand{\D}{\mathbb D}
\newcommand{\F}{\mathbf F}
\newcommand{\RR}{\mathbb R}
\newcommand{\PP}{\mathbb P}

\newcommand{\Z}{\mathbb Z}
\newcommand{\I}{\mathbb I}
\newcommand{\A}{\mathbb A}

\begin{document}
\maketitle \small

\noindent \textbf{Abstract.} A linear system of differential
equations describing a joint motion of elastic porous body and
fluid occupying porous space is considered.  Although the problem
is linear, it is very hard to tackle due to the fact that its main
differential equations involve non-smooth oscillatory
coefficients, both big and small, under the differentiation
operators. The rigorous justification, under various conditions
imposed on physical parameters, is fulfilled for homogenization
procedures as the  dimensionless size of the pores tends to zero,
while the porous body is geometrically periodic. As the results,
we derive Biot's equations of poroelasticity, equations of
viscoelasticity, or decoupled  system consisting of non-isotropic
Lam\'{e}'s equations and Darcy's system of filtration, depending
on ratios between physical parameters. The proofs are based on
Nguetseng's two-scale
convergence method of homogenization in periodic structures.\\


\normalsize

\addtocounter{section}{1} \setcounter{equation}{0}

\begin{center} \textbf{Introduction}
\end{center}

In this article a problem of modelling of small perturbations in
elastic deformable medium, perforated by a system of channels
(pores) filled with liquid or gas, is considered. Such media are
called \textit{elastic porous media} and they are a rather good
approximation to real consolidated grounds. In present-day
literature, the field of study in mechanics corresponding to these
media is called \textit{poromechanics} \cite{OC}. The solid
component of such a medium has a name of \textit{ skeleton}, and
the domain, which is filled with a fluid, is named a
\textit{porous space}. The exact mathematical model of elastic
porous medium consists of the classical equations of momentum and
mass balance, which are stated in Euler variables, of the
equations determining stress fields in both solid and liquid
phases, and of an endowing relation determining behavior of the
interface between liquid and solid components. The latter relation
expresses the fact that the interface is a material surface, which
amounts to the condition that it consists of the same material
particles all the time. Denoting by $\rho$ the density of medium,
by ${\mathbf v}$ the velocity, by $\PP^f$ the stress tensor in the
liquid component, by $\PP^s$ the stress tensor in the rigid
skeleton, and by $\tilde{\chi}$ the characteristic (indicator)
function of porous space, we write the fundamental differential
equations of the nonlinear model in the form
\begin{equation} \nonumber
\rho \frac{d {\mathbf v}}{dt} =\div_x \{\tilde{\chi} \PP^{f}+
(1-\tilde{\chi} )\PP^{s}\}+ \rho \F,
\end{equation}
\begin{equation}\nonumber
   \frac{d\rho}{dt}+ \rho \div_x {\mathbf v} =0, \quad
   \frac{d\tilde{\chi}}{dt} =0,
\end{equation}
where $d/dt$ stands for the material derivative with respect to
the time variable.

Clearly the above stated original model is a model with an unknown
(free) boundary. The  more precise formulation of the nonlinear
problem is not in focus of our present work. Instead, we aim to
study the problem, linearized at the rest state. In continuum
mechanics the methods of linearization are developed rather
deeply. The so obtained linear model is a commonly accepted and
basic one for description of filtration and seismic acoustics in
elastic porous media (see, for example, \cite{B-K,S-P,G-M1}).
Further we refer to this model as to \textbf{model A}. In this
model the characteristic function of the porous space
$\tilde{\chi}$ is a known function for $t>0$. It is assumed that
this function coincides with the characteristic function of the
porous space $\bar{\chi}$, given at the initial moment. Being
written in terms of dimensionless variables, the differential
equations of the model involve frequently oscillating non-smooth
coefficients, which have structures of linear combinations of the
function $\bar{\chi}$. These coefficients undergo differentiation
with respect to $\x$ and besides may be very big or very small
quantities as compared to the main small parameter $\varepsilon$.
In the model under consideration we define $\varepsilon$ as the
characteristic size of pores $l$ divided by the characteristic
size $L$ of the entire porous body:
$$\varepsilon =\frac{l}{L}.$$

Denoting by ${\mathbf w}$ the dimensionless displacement vector of
the continuum medium, in terms of dimensionless variables we write
the differential equations of model A as follows:
\begin{equation} \nonumber
\alpha_\tau \bar{\rho}\frac{\partial ^{2}\w}{\partial t^{2}}
=\div_x \PP +\bar{\rho}\F,
\end{equation}
\begin{equation} \nonumber
\PP = \bar{\chi} \alpha_\mu \D(x,\frac{\partial \w}{\partial t})
+(1-\bar{\chi})\alpha_\lambda \D(x,\w) -(q+\pi)\I,
\end{equation}
\begin{equation*}
p=-\alpha_p \bar{\chi}\div_x \w,\quad \pi=-\alpha_\eta
(1-\bar{\chi})\div_x \w,\quad
q=p+\frac{\alpha_{\nu}}{\alpha_{p}}\frac{\partial p}{\partial t}.
\end{equation*}
Here and further we use the notation
\begin{equation*}
\bar{\rho}=\bar{\chi}\rho_f +(1-\bar{\chi})\rho_s, \quad
\D(x,\uu)= (1/2)\left(\nabla_x \uu +(\nabla_x \uu)^T\right).
\end{equation*}

From the purely mathematical point of view, the corresponding
initial-boundary value problem for model\textbf{ A }is well-posed
in the sense that it has a unique solution belonging to a suitable
functional space on any finite temporal interval. However, in view
of possible applications, for example, for developing numerical
codes, this model is ineffective due to its sophistication even if
a modern supercomputer is available. Therefore a question of
finding an effective approximate models is vital. Since the model
involves the small parameter $\varepsilon$, the most natural
approach to this question is to derive models that would describe
limiting regimes arising as $\varepsilon$ tends to zero. Such an
approximation significantly simplifies the original problem and at
the same time preserves all of its main features. But even this
approach is too hard to work out, and some additional simplifying
assumptions are necessary. In terms of geometrical properties of
the medium, the most appropriate is to simplify the problem
postulating that the porous structure is periodic. Further by
\textbf{model} ${\mathbf B}^\varepsilon$ we will call model
\textbf{A} supplemented by this periodicity condition. Thus, our
main goal now is a derivation of all possible homogenized
equations in the model ${\mathbf B}^\varepsilon$.

The first research with the aim of finding limiting regimes in the
case when the skeleton was assumed to be an absolutely rigid body
was carried out by E. Sanchez-Palencia and L. Tartar.  E.
Sanchez-Palencia \cite[Sec. 7.2]{S-P} formally obtained Darcy's
law of filtration  using the method of  two-scale asymptotic
expansions, and L. Tartar \cite[Appendix]{S-P} mathematically
rigorously justified the homogenization procedure. Using the same
method of  two-scale expansions J. Keller and R. Burridge
\cite{B-K} derived formally the system of Biot's equations
\cite{BIOT} from model ${\mathbf B}^\varepsilon$ in the case when
the   parameter $\alpha _{\mu}$  was of order $\varepsilon^2$, and
the rest of the coefficients were fixed independent of
$\varepsilon$. It is well-known that the various modifications of
Biot's model are bases of seismic acoustics problems up-to-date.
This fact emphasizes importance of comprehensive study of model
\textbf{A} and model ${\mathbf B}^\varepsilon$ one more time. J.
Keller and R. Burridge also considered model
$\mathbf{B}^\varepsilon$ under assumption that all the physical
parameters were fixed independent of $\varepsilon$, and formally
derived as the result a system of equations of viscoelasticity.

 Under the same assumptions as in the article \cite{B-K}, the rigorous
justification of Biot's model was given by G. Nguetseng \cite{GNG}
and later by A. Mikeli\'{c}, R. P. Gilbert, Th. Clopeaut, and J.
L. Ferrin in \cite{G-M1,G-M2,G-M3}.  Also A. Mikeli\'{c}
\textit{et al} derived  a system of equations of viscoelasticity,
when all the physical parameters were fixed independent of
$\varepsilon$. In these works, Nguetseng's two-scale convergence
method \cite{NGU,LNW} was the main method of investigation of the
model $\mathbf{B}^\varepsilon$.

In the present work by means of the same  method we investigate
all possible limiting regimes in the model
$\mathbf{B}^\varepsilon$. This method in rather simple form
discovers the structure of the weak limit of a sequence
$\{z^\varepsilon\}$ as $\varepsilon \searrow 0$, where
$z^\varepsilon =u^\varepsilon v^\varepsilon$ and sequences
$\{u^\varepsilon\}$ and $\{v^\varepsilon\}$ converge as
$\varepsilon \searrow 0$ merely weakly, but at the same time
function $u^\varepsilon$ has the special structure
$u^\varepsilon(\x) =u(\x/\varepsilon)$ with $u(\y)$ being periodic
in $\y$.

Moreover,  Nguetseng's  method allows to establish asymptotic
expansions of a solution of model $B^\varepsilon$ in the form
\begin{equation} \nonumber
{\w^{\varepsilon}}({\x},t)=\varepsilon^{\beta}\left({\w}_0({\x},t)+\varepsilon
\w_1\left(\x,t, \x / \varepsilon\right)+ \textit{o}(\varepsilon
)\right),
\end{equation}
where $\w_0(\x,t)$ is a solution of the homogenized (limiting)
problem, $\w_1 (\x,t,\y)$ is a solution of some initial-boundary
value problem posed on the generic periodic cell of the porous
space, and exponent $\beta$ is defined by dimensionless parameters
of the model. Distinct asymptotic behavior of these parameters and
distinct geometry of the porous space lead to different limiting
regimes, namely, to various forms of Darcy's law for velocity of
liquid component and of non-isotropic Lam\'{e}'s equations for
displacement of rigid component in cases of big parameter
$\alpha_\lambda$, also, to various forms of Biot's system in cases
of small parameter $\alpha_\mu$, and to different forms of
equations of viscoelasticity in cases when parameters $\alpha_\mu$
and $\alpha_\lambda$ are $O(1)$ as $\varepsilon \searrow 0$. For
example, in the case when
$$
\alpha_\lambda=\varepsilon^{-\frac{3}{2}}, \quad
\alpha_\mu=\varepsilon^{\frac{1}{2}},
$$
the velocity $\vv^{\varepsilon}=\chi^{\varepsilon}
\partial \w^\varepsilon / \partial t$ of liquid component and the
displacement
$\uu^{\varepsilon}=(1-\chi^{\varepsilon})\w^{\varepsilon}$ of
rigid skeleton possess the following asymptotic:
\begin{equation} \nonumber
\vv^\varepsilon(\x,t)=\varepsilon^{\frac{3}{2}}\chi^{\varepsilon}\left(\V\left(\x,t,
\x / \varepsilon\right)+ \textit{o}(\varepsilon)\right),
\end{equation}
\begin{equation} \nonumber
\uu^\varepsilon(\x,t)=
\varepsilon^{\frac{3}{2}}(1-\chi^{\varepsilon})\left(\uu({\x},t)+\varepsilon
\UU\left(\x,t, \x / \varepsilon\right)+ \textit{o}(\varepsilon
)\right).
\end{equation}
At the same time all equations are determined in a unique way by
the given physical parameters of the original model and by
geometry of the porous space. For example, in the case of isolated
pores (disconnected porous space) the unique limiting regime  for
any combinations of parameters is a regime described by the
non-isotropic system of Lam\'{e}'s equations.

On our opinion, the proposed approach, when the limiting
transition in all coefficients is fulfilled simultaneously, is the
most natural one. We emphasize that it is not assumed starting
from the original model that the fluid component is inviscid, the
porous skeleton is absolutely rigid, or either of the components
is incompressible. These kinds of properties arise in limiting
models depending on limiting relations, which involve all
parameters of the problem.

 The articles \cite{G-M1} and
\cite{G-M2}, as well as the present one, show in favor of such a
uniform approach, because they exhibit the situations, when
various rates of approach of the small parameter $\alpha_{\mu}$ to
zero yield distinct homogenized equations. Moreover, these
equations differ from the  homogenized  equations derived as the
limit of model $\mathbf{B}^\varepsilon$ under assumption
$\alpha_\mu\equiv 0$, imposed even before homogenization
\cite{G-M3}.

Suppose that all dimensionless parameters of the
$\mathbf{B}^\varepsilon$ depend on the small parameter
$\varepsilon$ of the model and there exist limits (finite or
infinite)
\begin{equation} \nonumber
\lim_{\varepsilon\searrow 0} \alpha_\mu(\varepsilon) =\mu_0, \quad
\lim_{\varepsilon\searrow 0} \alpha_\lambda(\varepsilon)
=\lambda_0, \quad \lim_{\varepsilon\searrow 0}
\alpha_\tau(\varepsilon)=\tau_{0},
\end{equation}
We restrict our consideration by the cases when $\tau_0<\infty$
and either of the following situations has place.
\begin{eqnarray*}
{\bf (I)\phantom{xxxxxxxx}} & \mu_0=0, &  0< \lambda_0 <\infty,\\
{\bf (II)\phantom{xxxxxxxx}} &0\leq \mu_0<\infty, & \lambda_0=\infty,\\
{\bf (III)\phantom{xxxxxxxx}} & 0<\mu_0, \lambda_0 <\infty.\\
\end{eqnarray*}
If  $\tau_0=\infty$ then, re-normalizing the displacement vector
by setting
\begin{equation}\nonumber
\w  \rightarrow  \alpha_\tau \w,
\end{equation}
we reduce the problem to one of the cases {\it (I)--(III)}.

In the present paper we show that in the case ${\bf (I)}$ the
homogenized equations have various forms of Biot's system of
equations of poroelasticity for a two-velocity continuum  media or
non-isotropic Lam'{e}'s system of equations of one-velocity
continuum  media (for example, for the case of disconnected porous
space)(theorem \ref{theorem2.2}). In the case ${\bf (II)}$ the
homogenized equations are different modifications of Darcy's
system of equations of filtration for the velocity of the liquid
component (and as a first approximation the solid component
behaves yourself as an absolutely rigid body) and as a second
approximation -- non-isotropic Lam'{e}'s system of equations for
the re-normalized displacements of the solid component or Biot's
system of equations of poroelasticity for the re-normalized
displacements of the liquid and solid components(theorem
\ref{theorem2.3}). Finally, in the case ${\bf (III)}$ they are the
non-local viscoelasticity equations or non-isotropic non-local
Lam'{e}'s system of equations of one-velocity continuum media
(theorem \ref{theorem2.4}).

 \setcounter{equation}{0}
\setcounter{theorem}{0} \setcounter{lemma}{0}
\setcounter{proposition}{0} \setcounter{corollary}{0}
\setcounter{definition}{0} \setcounter{assumption}{0}
\begin{center} \textbf{\S1. Models   A and} \mathbf{B}$^\varepsilon$
\end{center}

\textbf{1.1. Differential equations, boundary and initial
conditions.} Let a domain $\Omega^*$ of the physical space $\RR^3$
be the union of a domain $\Omega_s^*$ occupied with the rigid
porous ground and a domain $\Omega_f^*$ corresponding to hollows
(pores) in the ground. Domain $\Omega_f^*$ is called a
\textit{porous space} and is assumed to be filled with liquid or
gas. Denote by ${\mathbf w}^*({\mathbf x}^*,t^*)$ the displacement
vector of the continuum medium (of ground or liquid or gas)) at
the point ${\mathbf x}^*\in \Omega^*$ in Euler coordinate system
at the moment of time $t^*\geq 0$. Under an assumption that the
displacement vector ${\mathbf w}^*$ is small in $\Omega^*$, which
amounts to a case when deformations  are small, dynamics of rigid
phase is described by linear Lam\'{e}'s equations and dynamics of
fluid or gas is described by Stokes equations. At the same time we
may set that the velocity vector ${\mathbf v}^*_f({\mathbf
x^*},t^*)$ in fluid or gas is a partial derivative of the
displacement vector with respect to time variable, i.e., that
${\mathbf v}^*_f({\mathbf x^*},t^*)=\partial {\mathbf
w}^*_f({\mathbf x}^*,t^*)/\partial t^*$. This assumption makes
perfect sense in description of continuous media in domains, where
the characteristic size of pores $l$ is very small as compared to
the diameter $L$ of the domain $\Omega^*$, i.e., $l=\varepsilon L$
and $\varepsilon<<1$ (see, for example, \cite{B-K,S-P,G-M1} and
recall the observation of the linearization procedure for the
exact nonlinear model in Introduction).

In terms of the dimensionless variables, not denoted by the
asterisk  $_*$  below,
\begin{equation} \label{(DL1)}
\left. \begin{array}{r} \displaystyle {\mathbf x}^* = {\mathbf x}
L,\quad t^* = t \tau, \quad {\mathbf w}^*_f = {\mathbf w}_f
L,\quad {\mathbf w}^*_s = {\mathbf w}_s L, \\[1ex]
\displaystyle {\mathbf F}^* = {\mathbf F} g, \quad p^*_f = p
p_0,\quad \rho^*_f = \rho_f \rho_0, \quad \rho^*_s = \rho_s
\rho_0,
\end{array} \right\}
\end{equation}
the displacement ${\mathbf w}_f({\mathbf x},t)$ and the pressure
$p({\mathbf x},t)$ of fluid and the displacement ${\mathbf
w}_s({\mathbf x},t)$ of rigid skeleton satisfy the system of
Stokes equations
\begin{eqnarray} \label{(S1)}
\displaystyle \alpha_\tau \rho_f\frac{\partial^2 {\mathbf
w}_f}{\partial t^2}
 =\div_x \PP^f + \rho_f \F, & & \x\in \Omega_f,\; t>0, \\
\label{(S2)}  \displaystyle  \PP^f =- p \I + \alpha_\nu
\left(\div_x \frac{\partial \w_f}{\partial t}\right)\cdot \I+
\alpha_\mu  \D\left(x,\frac{\partial \w_f}{\partial t}\right),
& & \x\in \Omega_f,\; t>0,\\
\label{(3)} \displaystyle p =-\alpha_p  \div_x \w_f, & & \x\in
\Omega_f,\; t>0
\end{eqnarray}
and the system of Lam\'{e}'s equations
\begin{eqnarray} \label{(L1)}
\displaystyle \alpha_\tau \rho_s \frac{\partial^2 \w_s}{\partial
t^2} =\div_x \PP^s +  \rho_s \F, & & \x\in \Omega_s,\; t>0,\\
\label{(L2)} \displaystyle \PP^s = \alpha_\eta  (\div_x \w_s)\cdot
\I + \alpha_\lambda  \D(x,\w_s), & & \x\in \Omega_s,\; t>0.
\end{eqnarray}
In  \eqref{(DL1)}--\eqref{(L2)} $\F=\F(\x,t)$ is the given vector
of distributed mass forces, $L$ is the characteristic macroscopic
size -- the diameter of the domain $\Omega^*$, $\tau$ is
characteristic duration of physical processes, $\rho_0$ is mean
density of air under atmosphere pressure, $\rho_f$ and $\rho_s$
are respectively mean dimensionless densities of liquid and rigid
phases, correlated with mean density of air, $g$ is the value of
acceleration of gravity, and $p_0$ is atmosphere pressure.

Dimensionless constants $\alpha_i$ $(i=\tau,\nu,\ldots)$ are
defined by the formulas
\begin{equation} \label{(DL2)}
\left. \begin{array}{r} \displaystyle \alpha_\tau =\frac{ L}{g
\tau^2}, \quad \alpha_\nu = \frac{\nu }{\tau Lg\rho_0}, \quad
\alpha_\mu =\frac{2\mu }{\tau Lg\rho_0},\\[1ex] \displaystyle
 \quad \alpha_p =\frac{c^{2}\rho _{f}}{ Lg},
 \quad \alpha_\eta =\frac{\eta}{ Lg\rho_0},
 \quad \alpha_\lambda =\frac{2 \lambda }{ Lg\rho_0},
\end{array} \right\}
\end{equation}
where $\mu$ is the  viscosity of fluid or gas, $\nu$ is the bulk
viscosity of fluid or gas, $\lambda$ and $\eta$ are elastic
Lam\'{e}'s constants, and  $c$ is a  speed of sound in fluid.

For the unknown functions $\w_f$, $p$, and $\w_s$, the commonly
accepted conditions of continuity of the displacement field and
normal tensions are imposed on the interface $\Gamma:=\partial
\Omega_f\cap \partial \Omega_s$ between the two phases (see, for
example, \cite{B-K,S-P,G-M1}):
\begin{eqnarray} \label{(RH1)}
\displaystyle \PP^f\cdot{\mathbf n} = \PP^s \cdot {\mathbf n}, & &
\x \in \Gamma,\; t>0,\\
\label{(RH2)} \displaystyle \w^f =\w^s, & & \x \in \Gamma,\; t>0.
\end{eqnarray}
 In \eqref{(RH1)} ${\mathbf n}$ is the unit normal vector to $\Gamma$.
 Note that exactly these conditions appear as the result of linearization
 of the exact nonlinear model.

Finally, system \eqref{(S1)}--\eqref{(L2)},
\eqref{(RH1)}--\eqref{(RH2)} is endowed by giving a displacement
field on $S=\partial \Omega$ and at the moment $t=0$ and by giving
a velocity field at the $t=0$. Further without loss of generality,
with the aim to simplify the technical outline, we suppose that
these boundary conditions are homogeneous.

\textbf{1.2. Geometry of porous space.} In model \textbf{A }a
Lipschitz smoothness of the interface between porous space and
rigid skeleton is the only restriction on geometry of porous
space. In model $\mathbf{B}^\varepsilon$ the porous medium has
geometrically periodic structure. Its formal description is as
follows \cite{G-M1,B-P}.

Firstly a geometric structure inside a pattern unit cell
$Y=(0,1)^3$ is defined. Let $Y_s$ be a `solid part' of the cell
$Y$. The 'liquid part' $Y_f$ is its open complement. Set
$Y_s^k:=Y_s +{\mathbf k}$, ${\mathbf k}\in \Z^3$, the translation
of $Y_s$ on an integer-valued vector ${\mathbf k}$. Union of such
translations along all ${\mathbf k}$, $E_s:= \cup_{{\mathbf k}\in
\Z^3} Y_s^k$ is the 1-periodic repetition of $Y_s$ all over
$\RR^3$. Let $E_f$ be the open complement of $E_s$ in $\RR^3$. The
following assumptions on geometry of $Y_f$ and $E_f$ are
accepted.\\[1ex]
(i) $Y_s$ is an open connected set of strictly positive measure
with a Lipschitz boundary, and $Y_f$ also has strictly positive
measure on $Y$.\\[1ex]
(ii) $E_f$ and $E_s$ are open sets with $C^{0,1}$-smooth
boundaries. The set $E_f$ is locally situated on one side of the
boundary $\partial E_f$, and the \textbf{set} $\textbf{E}_s$ is
locally situated on one side of the boundary $\partial E_s$  and
\textbf{connected}.

Domains $\Omega^\varepsilon_s$ and $\Omega_f^\varepsilon$ are
intersections of the domain $\Omega$ with the sets $\RR^3_s$ and
$\RR^3_f$, where the sets $\RR^3_s$ and $\RR^3_f$ are periodic
domains in $\RR^3$ with generic cells $\varepsilon Y_s$ and
$\varepsilon Y_f$ of the diameter $\varepsilon$, respectively.

Union $\bar{Y}_s \cup \bar{Y}_f$ is the closed cube
$\bar{Y}=\{\y\in \RR^3,\; 0\leq y_i\leq 1,\; i=1,2,3\}$, and the
interface $\Gamma^\varepsilon =\partial \Omega_s^\varepsilon \cap
\partial \Omega_f^\varepsilon$ is the $\varepsilon$-periodic repetition of
the boundary $\varepsilon \gamma =\varepsilon \partial Y_f \cap
\partial Y_s$ all over $\omega$.

Further by $\bar{\chi}=\chi^\varepsilon$ we will denote the
characteristic function of the porous space.

For simplicity we accept the following constraint on the domain
$\Omega$ and the parameter $\varepsilon$.
\begin{assumption} \label{assumption1}
Domain $\Omega$ is cube, $\Omega:=(0,1)^3$, and quantity
$1/\varepsilon$ is integer, so that $\Omega$ always contains an
integer number of elementary cells $Y_i^\varepsilon$.
\end{assumption}
Under this assumption, we have
\begin{equation} \label{char}
 \bar{\chi}(\x)=\chi^{\varepsilon}(\x)=\chi
 \left(\x / \varepsilon\right),
\end{equation}
where $\chi (\y)$ is the characteristic function of $Y_f$ in $Y$.

 We say that a \textbf{porous space is disconnected (isolated
 pores)}
 if  $\gamma \cap \partial Y=\emptyset$.\\

\textbf{1.3. Generalized solutions in models A and}
$\mathbf{B}^\varepsilon $. Define the displacement $\w(\x,t)$ in
the whole domain $\Omega$ by the formula
\begin{equation} \label{(12)}
\w(\x,t)=\left\{\begin{array}{l}  \w_f(\x,t),\quad \x\in
\Omega_f,\quad t>0,\\  \w_s(\x,t),\quad \x\in \Omega_s,\quad t>0,
\end{array} \right.
\end{equation}
and the pressures $p(\x,t)$, $q(\x,t)$, and $\pi(\x,t)$ by
formulas
\begin{equation} \label{(13)}
p=-\alpha_p \bar{\chi}\div_x \w,\quad \pi=-\alpha_\eta
(1-\bar{\chi})\div_x \w,\quad
q=p+\frac{\alpha_{\nu}}{\alpha_{p}}\frac{\partial p}{\partial t}.
\end{equation}
Thus introduced new unknown functions should satisfy the system
\begin{subequations} \label{(14)}
\begin{equation} \label{(14a)}
\alpha_\tau \bar{\rho}\frac{\partial ^{2}\w}{\partial t^{2}}
=\div_x \PP +\bar{\rho}\F,
\end{equation}
\begin{equation} \label{(14b)}
\PP = \bar{\chi} \alpha_\mu \D(x,\frac{\partial \w}{\partial t})
+(1-\bar{\chi})\alpha_\lambda \D(x,\w) -(q +\pi )\I,
\end{equation}
\end{subequations}
 where $\bar{\rho}=\bar{\chi} \rho_f
+(1-\bar{\chi})\rho_s$.  If
$\bar{\chi}(\x)=\chi^{\varepsilon}(\x)$, then $\bar{\rho}=\rho
^{\varepsilon} \equiv \chi^{\varepsilon} \rho_f
+(1-\chi^{\varepsilon})\rho_s$.

Equations \eqref{(14)} are understood in the sense of
distributions theory. They involve the both equations \eqref{(S1)}
and \eqref{(L1)} in the domains $\Omega_f$ and $\Omega_s$,
respectively, and the boundary conditions \eqref{(RH1)} and
\eqref{(RH2)} on the interface $\Gamma$. There are various
equivalent in the sense of distributions forms of representation
of equations \eqref{(14)}. In what follows, it is convenient to
write them in the form of the integral equality
\begin{eqnarray} \label{(15)}
&&\int_{\Omega_{T}} \{\alpha_\tau \bar{\rho} \w \cdot
\frac{\partial ^{2}{\mathbf \psi}}{\partial t^{2}} -\bar{\rho} \F
\cdot {\mathbf \psi} -\bar{\chi}
\alpha_\mu \D(x,\w):\D(x,\frac{\partial {\mathbf \psi}}{\partial t})+\\
&&((1-\bar{\chi})\alpha_\lambda \D(x,\w) -(q +\pi
)\I):\D(x,{\mathbf \psi})\} d\x dt =0,\nonumber
\end{eqnarray}
where ${\mathbf \psi}(\x,t)$ is an arbitrary smooth test function,
such that ${\mathbf \psi}=\partial {\mathbf \psi}/ \partial t=0$
at the  $t=T$ and ${\mathbf \psi}=0$  on the boundary $S_T$ of the
domain $\Omega_T$, $\Omega_T=\Omega \times (0,T)$.

In \eqref{(15)} by $A:B$ we denote the convolution (or,
equivalently, the inner tensor product) of two second-rank tensors
along the both indexes, i.e., $A:B=\mbox{tr\,} (B^*\circ
A)=\sum_{i,j=1}^3 A_{ij} B_{ji}$.

Mathematical model \eqref{(13)}, \eqref{(14b)}, and \eqref{(15)}
is endowed by the initial and boundary conditions
\begin{equation} \label{(16)}
\w(\x,0)=0, \quad \frac{\partial \w}{\partial t}(\x,0)=0,\quad
\x\in \Omega,
\end{equation}
\begin{equation} \label{(17)}
 \w(\x,t)=0, \quad \x\in S, \quad t>0.
\end{equation}
Here the assumption, that the boundary and initial conditions are
homogeneous, is not essential.

We will call solution of problem \eqref{(13)}, \eqref{(14b)},
\eqref{(15)}--\eqref{(17)} a \textit{generalized solution in
\textbf{model A}} and, if $\Omega$ is supplemented with periodic
structure, a \textit{generalized solution in \textbf{model}
$\textbf{B}^\varepsilon$}.

The following assertion states that problem \eqref{(13)},
\eqref{(14b)}, \eqref{(15)}--\eqref{(17)} is well-posed.
\begin{lemma} \label{Lemma1.1}
Let the interface $\Gamma$ between $\Omega_s$ and $\Omega_f$ be
piece-wise continuously differentiable, parameters $\rho_f$,
$\rho_s$, $\alpha_{\tau}$, $\alpha_\mu$, $\alpha_\nu$,
$\alpha_\lambda$, $\alpha_\eta$, and $\alpha_p$ be strictly
positive, and assume that $\F,
\partial \F / \partial t \in L^2(\Omega_T)$.

Then, on any finite temporal interval $[0,T]$, problem
\eqref{(13)}, \eqref{(14b)}, \eqref{(15)}--\eqref{(17)} has a
unique solution $\{\w,p,\pi \}$, and this solution admits the
bounds
\begin{equation*}
\max\limits_{0<t<T}(\sqrt{\alpha_\eta}\|(1- \bar{\chi}) \div_x
 \frac{\partial \w}{\partial t}(t) \|_{2,\Omega}+
  \sqrt{\alpha_\lambda}\|(1-\bar{\chi})
\nabla_x \frac{\partial \w}{\partial t}(t) \|_{2,\Omega}
\end{equation*}
\begin{equation*}
+ \sqrt{\alpha_\tau}\|\frac{\partial ^2\w}{\partial
t^2}(t)\|_{2,\Omega}+ \sqrt{\alpha _{p}}\| \bar{\chi} \div_x
\frac{\partial \w}{\partial t}(t)\|_{2,\Omega })+ \sqrt{\alpha
_{\nu}}\| \bar{\chi} \div_x\frac{\partial ^2\w}{\partial
t^2}\|_{2,\Omega _{T}}+
\end{equation*}
\begin{equation} \label{(18)}
+\sqrt{\alpha_\mu}\|\bar{\chi} \nabla_x\frac{\partial
^2\w}{\partial t ^2} \|_{2,\Omega_T} \leq
\frac{C}{\sqrt{\alpha_\tau}}\||\F|+|\frac{\partial \F}{\partial
t}|\|_{2,\Omega_{T}},
\end{equation}
where $C$ depends only on $T$.
\end{lemma}

Due to linearity of the problem, justification of Lemma
\ref{Lemma1.1} reduces to verification of bounds \eqref{(18)}.
These appear by means of differentiation of Eqs. \eqref{(14)} with
respect to $t$ (note that $\chi$ and $\bar{\rho}$ do not depend on
$t$), multiplication of the resulting equation by $\partial ^2\w /
\partial t^2$, and integration by parts. Pressures $q$ and $\pi$
are estimated directly from Eqs. \eqref{(13)}.

Further the focus of this article is centered solely on model
$\mathbf{B}^\varepsilon$, in which coefficients of Eqs.
\eqref{(13)} and \eqref{(14)} depend continuously on the small
parameter $\varepsilon$, $\bar{\rho}=\rho
^{\varepsilon}(x)=\chi^{\varepsilon}\rho_{f}+(1-\chi^{\varepsilon})\rho_{s}$
and  $\{\w^{\varepsilon}, p^{\varepsilon},
q^{\varepsilon},\pi^{\varepsilon} \}$ is a corresponding
generalized solution. We aim to find out the limiting regimes of
the model as $\varepsilon \searrow 0$.

\addtocounter{section}{1}\setcounter{equation}{0}
\setcounter{theorem}{0} \setcounter{lemma}{0}
\setcounter{proposition}{0} \setcounter{corollary}{0}
\setcounter{definition}{0} \setcounter{assumption}{0}
\begin{center} \textbf{\S2. Formulation of the main results}
\end{center}
Suppose additionally that there exist limits (finite or infinite)
\begin{equation*}
\lim_{\varepsilon\searrow 0} \alpha_\nu(\varepsilon) =\nu_0, \quad
\lim_{\varepsilon\searrow 0} \alpha_p(\varepsilon) =p_{*}, \quad
\lim_{\varepsilon\searrow 0} \alpha_\eta(\varepsilon) =\eta_0.
\end{equation*}
\begin{equation} \nonumber
\lim_{\varepsilon\searrow 0} \frac{\alpha_\mu}{\varepsilon^{2}}
=\mu_1, \quad \lim_{\varepsilon\searrow 0} \frac{\varepsilon^{2}
\alpha_p}{\alpha_\mu}=p_{1}, \quad \lim_{\varepsilon\searrow 0}
\frac{\alpha_\lambda \varepsilon^{2}}{\alpha_\mu} =\lambda_{1},
\end{equation}
\begin{equation} \nonumber
 \lim_{\varepsilon\searrow 0}
\frac{\alpha_\eta \varepsilon^{2}}{\alpha_\mu}=\eta_{1},
 \quad ,\lim_{\varepsilon\searrow 0} \frac{\alpha_\eta}{\alpha_\lambda} =\eta_{2},
\quad \lim_{\varepsilon\searrow 0} \frac{\alpha_p}{\alpha_\lambda}
=p_{2}.
\end{equation}
 In what follows we assume that
\begin{assumption} \label{assumption3}
Dimensionless parameters in the model $\mathbf{B}^{\varepsilon}$
satisfy restrictions
\begin{equation}\label{2.01}
  p_{*}^{-1}, \quad \mu_{0}, \quad \nu _{0}, \quad \lambda_{0}^{-1} <\infty;
   \quad   0< \tau _{0}+ \mu_1.
\end{equation}
\end{assumption}

 All parameters may take all permitted values. For example, if
 $\tau_{0}=0$ or  $p_{*}^{-1}=0$, then all terms in final equations
 containing these  parameters  disappear.

 The following Theorems
\ref{theorem2.1}--\ref{theorem2.4} are the main results of the
paper.

\begin{theorem} \label{theorem2.1}
Assume that conditions of Lemma \ref{Lemma1.1} hold and that \\
 $\{\w^{\varepsilon},q^{\varepsilon}, p^{\varepsilon},\pi^{\varepsilon} \}$
is a generalized solution in model $\mathbf{B}^\varepsilon$.

The following assertions are valid:
\\[1ex]
\textbf{(I)} If
$$ \lambda _{0} <\infty, $$
 then
\begin{equation} \label{(2.2)}
 \displaystyle \max\limits_{0\leq t\leq
T}\| |\w^{\varepsilon}(t)|+ \sqrt{\alpha_\mu} \chi^\varepsilon
|\nabla_x \w^{\varepsilon}(t)|+  (1-\chi^\varepsilon) |\nabla_x
 \w^{\varepsilon}(t)| \|_{2,\Omega} \leq I_{F} ,
\end{equation}
\begin{equation}\label{(2.4)}
\|q^{\varepsilon}\|_{2,\Omega_{T}} +
\|p^{\varepsilon}\|_{2,\Omega_{T}} + \|\pi
^{\varepsilon}\|_{2,\Omega_{T}} \leq I_{F},
\end{equation}
where  $ I_{F}=C \||\F|+|\partial \F / \partial
t|\|_{2,\Omega_{T}}$ and  $C$ is a
constant  independent of $\varepsilon$.\\[1ex]
\textbf{(II)}   If
\begin{equation*}
\lambda_{0}=\infty,  \quad \mu_{1}=\infty,  \quad  0<\lambda_{1}
<\infty ,
\end{equation*}
then estimates  \eqref{(2.2)},  \eqref{(2.4)} hold for
re-normalized displacements
\begin{equation*}
    w^{\varepsilon} \rightarrow  \varepsilon^{-2}\alpha_{\mu} w^{\varepsilon},
\end{equation*}
with re-normalized parameters
 $$\alpha_{\mu} \rightarrow \varepsilon^{2},  \quad
 \alpha_{\lambda} \rightarrow \varepsilon^{2}\frac{\alpha_{\lambda }}{\alpha_{\mu}},
 \quad \alpha_{\tau} \rightarrow \varepsilon^{2}\frac{\alpha_{\tau }
 }{\alpha_{\mu}}, \quad \alpha_{\nu} \rightarrow
\varepsilon^{2}\frac{\alpha_{\nu}}{\alpha_{\mu}},\quad \alpha_{ p}
\rightarrow \varepsilon^{2}\frac{\alpha_{p}}{\alpha_{\mu}}.$$
\\[1ex]
\textbf{(III)}  If
\begin{equation*}
 \lambda_{0}=\infty,  \quad \mu_{1}<\infty,
\end{equation*}
then for displacements $\w^{\varepsilon}$ hold true estimates
 \eqref{(2.2)}  and under  condition
\begin{equation}\label{(2.6)}
  p_{*}<\infty,
\end{equation}
for the pressures  $q^{\varepsilon}$  and  $p^{\varepsilon}$ in
the liquid component hold true estimates  \eqref{(2.4)}.

 If instead of restriction \eqref{(2.6)} hold true conditions
\begin{equation}\label{(2.9.1)}
0< p_{2}; \quad\F=\nabla \Phi , \quad \frac{\partial
\Phi}{\partial t}, \quad  |\frac{\partial \F}{\partial t}|  \in
L^2(\Omega_T) ,
\end{equation}
then
\begin{equation}\label{(2.9.2)}
\displaystyle  \max\limits_{0\leq t\leq T}(\|(1-\chi^\varepsilon)
\nabla_x (\alpha_\lambda \w^{\varepsilon}(t)) \|_{2,\Omega}+\|
\chi^\varepsilon \div_x
(\alpha_\lambda\w^{\varepsilon}(t))\|_{2,\Omega }
   \leq I^{(1)}_{F},
\end{equation}
where  $I^{(1)}_{F}=C \||\F|+|\partial \Phi / \partial
t|+|\partial \F / \partial t|\|_{2,\Omega_{T}}$  and $C$ is a
constant  independent of $\varepsilon$.

These last  estimates   imply \eqref{(2.4)}.
\end{theorem}

Note, that for the last case $ \{\lambda_{0}=\infty,
 \mu_{1}<\infty \}$ we can get same estimates
\eqref{(2.2)}, \eqref{(2.4)}  and \eqref{(2.9.2)},  if instead
restrictions \eqref{(2.9.1)}  we consider
\begin{assumption} \label{assumption2}
$$\F=\F^{\varepsilon}(\x,t)(1- \chi ^{\varepsilon})$$
  and sequences
$\{\F^{\varepsilon}\}$  and $\{\partial\F^{\varepsilon} / \partial
t\}$ are uniformly bounded with respect to $\varepsilon$ in
$L^2(\Omega _{T})$.
\end{assumption}

\begin{theorem} \label{theorem2.2}
Assume that the hypotheses in Theorem \ref{theorem2.1} hold, and
\begin{equation*}
 \lambda_{0} <\infty, \quad  \mu_{0}=0.
\end{equation*}
Then functions $\w^{\varepsilon}$ admit an extension
 $\uu^{\varepsilon}$ from $\Omega_{s,T}^{\varepsilon}=\Omega_s^\varepsilon \times (0,T)$
 into $\Omega_{T}$
 such that the sequence $\{\uu^{\varepsilon}\}$ converges strongly
 in $L^{2}(\Omega_{T})$ and weakly in
 $L^{2}((0,T);W^1_2(\Omega))$ to the
 functions $\uu$. At the same time, sequences $\{\w^\varepsilon\}$,
 $\{p^{\varepsilon}\}$, $\{q^{\varepsilon}\}$, and $\{\pi^{\varepsilon}\}$ converge weakly
 in $L^{2}(\Omega_{T})$ to $\w$, $p$, $q$, and $\pi$, respectively.

 The following assertions for these limiting functions hold
 true:\\[1ex]
\textbf{(I)} If $\mu_1=\infty$ or the porous space is disconnected
(a case of isolated pores), then  $\w=\uu$ and  the functions
$\uu$, $p$, $q$, and $\pi$ satisfy in the domain $\Omega_{T}$  the
following initial-boundary value problem:
 \begin{equation}\label{(2.10)}
\tau _{0}\hat{\rho}\frac{\partial ^2\uu}{\partial t^2}=
    \div_x \{\lambda _{0}\A^{s}_{0}:\D(x,\uu) + B^{s}_{0}\div_x \uu
 +B^{s}_{1}q - (q+\pi )\cdot \I \}+\hat{\rho}\F,
\end{equation}
\begin{equation}\label{(2.11)}
\frac{1}{\eta_{0}}\pi+C^{s}_{0}:\D(x,\uu)+ a^{s}_{0}\div_x \uu
 +a^{s}_{1}q=0,
\end{equation}
\begin{equation}\label{(2.11.1)}
 \frac{1}{p_{*}}p+ \frac{1}{\eta_{0}}\pi + \div_x \uu=0,\quad
 p +\nu_0 p_{*}^{-1} \frac{\partial p}{\partial t}= q,
\end{equation}

where $\hat{\rho}=m \rho_{f} + (1-m)\rho_{s}$,
$m=\langle\chi\rangle$. The symmetric strictly  positively defined
constant fourth-rank tensor $\A^{s}_{0}$, matrices  $C^{s}_{0},
B^{s}_{0}$ and  $B^{s}_{1}$ and constants $a^{s}_{0}$ and
$a^{s}_{1}$  are defined below by formulas \eqref{(5.26)},
\eqref{(5.27.1)} - \eqref{(5.27.3)}.

Differential equations \eqref{(2.10)}--\eqref{(2.11)} are endowed
with the homogeneous initial and boundary conditions
 \begin{equation}\label{(2.12)}
 \tau _{0}\uu|_{t=0}= \tau _{0}\frac{\partial \uu}{\partial t}|_{t=0}=0,
\quad \x\in \Omega,\quad \uu(\x,t)=0, \quad \x\in S.
\end{equation}

\noindent
 \textbf{(II)} If $\mu_{1}<\infty$, then the weak limits $\uu$,
$\w^{f}$, $p$, $q$, $\pi$ of sequences $\{\uu^\varepsilon\}$,
$\{\chi ^{\varepsilon}\w^\varepsilon\}$,
  $\{p^\varepsilon\}$,  $\{q^\varepsilon\}$,
   $\{\pi^\varepsilon\}$ satisfy the initial-boundary
value problem consisting of the balance of momentum equation
\begin{eqnarray}\label{(2.15)}
&&\tau _{0}\rho_{f}\frac{\partial \vv}{\partial t}+\tau
_{0}\rho_{s}(1-m)
\frac{\partial ^2\uu}{\partial t^2} + \nabla (q+\pi )-\hat{\rho}\F=\\
&&\div_x \{\lambda _{0}\A^{s}_{0}:\D(x,\uu) + B^{s}_{0}\div_x \uu
 +B^{s}_{1}q \},\nonumber
\end{eqnarray}
where $\vv=\partial \w^{f} / \partial t$  and $\A^{s}_{0}$,
$B^{s}_{0}$ and $B^{s}_{1}$ are the same as in \eqref{(2.10)}, the
continuity equations \eqref{(2.11)},  the continuity equation and
the state equation
\begin{equation} \label{(2.16)}
  \frac{1}{p_{*}}p+ \frac{1}{\eta_{0}}\pi +\div_x \w^{f}  = (m-1)\div_x \uu,
\quad p +\nu_0 p_{*}^{-1} \frac{\partial p}{\partial t}= q,
\end{equation}

and Darcy's law in the form
\begin{equation}\label{(2.21)}
\vv=m\frac{\partial \uu}{\partial t}+\int_{0}^{t}
B_{1}(\mu_1,t-\tau)\cdot (-\nabla_x
q+\rho_{f}\F-\tau_{0}\rho_{f}\frac{\partial ^2 \uu}{\partial \tau
^2})(\x,\tau )d\tau
\end{equation}
in the case $\tau_{0}>0$ and $\mu_{1}>0$, Darcy's law in the form
\begin{equation}\label{(2.22)}
\vv=m\frac{\partial \uu}{\partial t}+B_{2}(\mu_{1})\cdot(-\nabla
q+\rho_{f}\F),
\end{equation}
in the case $\tau_{0}=0$ and $\mu_{1}>0$, and, finally, Darcy's
law in the form
\begin{equation}\label{(2.23)}
\vv=B_{3}\cdot \frac{\partial \uu}{\partial t}+\frac{1}{\tau
_{0}\rho_{f}}(m\I-B_{3})\cdot\int_{0}^{t}(-\nabla q(\x,\tau
)+\rho_{f}\F(\x,\tau ))d\tau,
\end{equation}
in the case $\tau_{0}>0$ and $\mu_{1}=0$.

This problem is endowed with initial and boundary conditions
\eqref{(2.12)} and the boundary condition
\begin{equation}\label{(2.24)}
\vv(\x,t)\cdot n(\x)=0,
     \quad \x \in S, \quad t>0,
\end{equation}
for the velocity $\vv$ of the fluid component.

In \eqref{(2.21)}--\eqref{(2.24)} $n(\x)$ is the unit normal
vector to $S$ at a point $\x \in S$, and matrices
$B_{1}(\mu_{1},t), B_{2}(\mu_{1})$ and $(m\I-B_{3})$ are defined
below by formulas \eqref{(5.33)}, \eqref{(5.35)} and
\eqref{(5.38)}.
\end{theorem}

\begin{theorem} \label{theorem2.3}
Assume that the hypotheses in Theorem \ref{theorem2.1} hold, and
that
 \begin{equation*}
  \lambda_{0} =\infty.
\end{equation*}
$\textbf{(I)}$  If $\mu_{1}<\infty$ and  one of conditions
  \eqref{(2.6)} or \eqref{(2.9.1)} holds  true, then sequences
 $\{\chi ^{\varepsilon}\w^\varepsilon\}$,
 $\{p^{\varepsilon}\}$  and $\{q^{\varepsilon}\}$ converge weakly
  in  $L^{2}(\Omega_{T})$ to  $\w^{f}$, $p$,  and  $q$  respectively.
 The functions $\w^{\varepsilon}$ admit an extension
 $\uu^{\varepsilon}$ from $\Omega_{s}^{\varepsilon}\times (0,T)$
  into $\Omega_{T}$
 such that the sequence $\{\uu^{\varepsilon}\}$ converges strongly
 in $L^2(\Omega_{T})$ and weakly in
 $L^{2}((0,T);W^1_2(\Omega))$ to zero
 and

1) if  $\tau_{0}>0$ and  $\mu_{1}>0$, then functions $\vv=\partial
\w^{f} / \partial t$, $p$ and  $q$ solve  in the domain
$\Omega_{T}$ the problem  $(F_{1})$, where
\begin{equation}\label{(2.25)}
\vv=\int_{0}^{t}B_{1}(\mu_{1},t-\tau)\cdot (-\nabla q(\x,\tau
)+\rho_{f}\F(\x,\tau ))d\tau,
\end{equation}
\begin{equation} \label{(2.26)}
    p +\nu_0 p_{*}^{-1} \frac{\partial p}{\partial t}= q,
     \quad  \frac{1}{p_{*}}\frac{\partial p}{\partial t}+\div_x \vv =0;
\end{equation}

2) if  $\tau_{0}=0$  and  $\mu_{1}>0$, then  functions $\vv$, $p$
 and  $q$ solve  in the domain $\Omega_{T}$ the problem  $(F_{2})$, where  $\vv$
 satisfies Darcy's law in the form
\begin{equation}\label{(2.29)}
\vv= B_{2}(\mu_{1})\cdot(-\nabla q+\rho_{f}\F),
\end{equation}
 and pressures  $p$ and  $q$ satisfy equations
 \eqref{(2.26)};

 finally,

3) if  $\tau_{0}>0$ and  $\mu_{1}=0$, then functions $\vv$, $p$
 and   $q$ solve  in the domain $\Omega_{T}$ the problem  $(F_{3})$,  where  $\vv$
satisfies Darcy's law in the form
\begin{equation}\label{(2.30)}
\vv= \frac{1}{\tau
_{0}\rho_{f}}(m\I-B_{3})\cdot\int_{0}^{t}(-\nabla q(\x,\tau
)+\rho_{f}\F(\x,\tau ))d\tau,
\end{equation}
 and pressures  $p$ and  $q$ satisfy equations
 \eqref{(2.26)}.

Problems $F_1$--$F_3$ are endowed with boundary condition
\eqref{(2.24)}.

 $\textbf{(II)}$ If $\mu_{1}<\infty$ and conditions  \eqref{(2.9.1)} hold
 true,
 then the sequence  $\{\alpha_{\lambda}\uu^\varepsilon \}$
 converges strongly in  $L^2(\Omega_{T})$ and weakly in
$L^2((0,T);W^1_2(\Omega))$ to function $\uu$ and the sequence
$\{\pi^{\varepsilon}\}$   converges weakly in $L^{2}(\Omega_{T})$
to the function  $\pi$.  The limiting functions $\uu$ and $\pi$
satisfy  the boundary value problem in the domain $\Omega $

 \begin{equation}\label{(2.31)}
0=\div_x \{\A^{s}_{0}:\D(x,\uu) + B^{s}_{0}\div_x \uu
 +B^{s}_{1}q - (q+\pi )\cdot \I \}+\hat{\rho}\F,
\end{equation}
\begin{equation}\label{(2.32)}
  \frac{1}{\eta_{2}}\pi+C^{s}_{0}:\D(x,\uu)+ a^{s}_{0}\div_x \uu
 +a^{s}_{1}q=0,
\end{equation}
where the function $q$ is referred to as given. It is defined from
the corresponding of Problems $F_1$--$F_3$  (the choice of the
problem depends on  $\tau_{0}$ and $\mu_{1}$). The symmetric
strictly  positively defined constant fourth-rank tensor
$\A^{s}_{0}$, matrices  $C^{s}_{0}, B^{s}_{0}$ and  $B^{s}_{1}$
and constants $a^{s}_{0}$ and  $a^{s}_{1}$  are defined below by
formulas \eqref{(5.26)}, \eqref{(5.27.1)} - \eqref{(5.27.2)}, in
which we have $\eta_0=\eta_2$  and  $\lambda_0 =1$.

This problem is endowed with the homogeneous  boundary
conditions.\\[1ex]
$\textbf{(III)}$  If  $\mu_{1}=\infty$,  $p_{1}^{-1}, \eta
_{1}^{-1}<\infty $ and $0< \lambda_{1}<\infty$, then there exist
weak limits
  $\w^{f}$, $p$ and $\pi$  of the sequences
 $\{\alpha_{\mu}\varepsilon^{-2}\chi ^\varepsilon\w^\varepsilon\}$,
  $\{p^\varepsilon\}$ and $\{\pi^\varepsilon\}$
  and a strong limit $\uu$  of the sequence $\{\alpha_{\mu}\varepsilon^{-2}\uu^\varepsilon
  \}$ in  $L^2(\Omega_{T})$ , which satisfy in $\Omega_{T}$ the following
   initial boundary-value problem:
\begin{equation}\label{(2.34)}
\left. \begin{array}{r} \displaystyle \div_x \{\lambda
_{1}\A^{s}_{0}:\D(x,\uu) + B^{s}_{0}\div_x \uu
 +B^{s}_{1}p - (p+\pi )\cdot \I \}+\hat{\rho}\F=0,\\[1ex]
\displaystyle \frac{\partial \w^{f}}{\partial t}=\frac{\partial
\uu}{\partial t}+B_{2}(1)\cdot (-\nabla p + \rho_{f}\F),\\[1ex]
\displaystyle \frac{1}{p_{1}} p+ \frac{1}{\eta _{1}}\pi + \div_x
\w^{f} =(m-1)\div_x \uu,\\[1ex]
\displaystyle \frac{1}{\eta_{1}}\pi+C^{s}_{0}:\D(x,\uu)+
a^{s}_{0}\div_x \uu +a^{s}_{1}p=0.
\end{array} \right\}
\end{equation}
 Here the symmetric strictly  positively defined constant fourth-rank tensor
$\A^{s}_{0}$, matrices  $C^{s}_{0}, B^{s}_{0}$ and  $B^{s}_{1}$
and constants $a^{s}_{0}$ and $a^{s}_{1}$  are defined below by
formulas \eqref{(5.26)}, \eqref{(5.27.1)} - \eqref{(5.27.3)}, in
which we have $\eta_{0}=\eta_{1}$ and  $\lambda_{0}=\lambda_{1}$.

This problem is endowed with the homogeneous  initial and boundary
conditions.

$\textbf{(IV)}$  If  $\mu_{1}=\infty$ and  $ \lambda_{1}=\infty$,
then the corresponding problem for displacements
$\{\alpha_{\mu}\varepsilon^{-2}\w^\varepsilon\}$  has been
considered in parts $\textbf{(I)}$-$\textbf{(II)}$ of the present
theorem.

\end{theorem}

\begin{theorem} \label{theorem2.4}
Assume that the hypotheses in Theorem \ref{theorem2.1} hold, and
that
$$0<\mu_{0}, \quad \lambda_{0}<\infty.$$
  Then  weak limits $\w$, $p$, $q$ and $\pi$ of
 sequences $\{\w^\varepsilon\}$,  $\{p^{\varepsilon}\}$, $\{\pi^\varepsilon\}$
 and  $\{q^{\varepsilon}\}$ satisfy in $\Omega_{T}$ the following initial-boundary
value problem:
 \begin{equation} \label{(2.35)}
\left. \begin{array}{r} \displaystyle \tau_0 \hat{\rho}
\frac{\partial ^2\w}{\partial t^2} + \nabla (q+\pi)-\hat{\rho}\F=\\[1ex]
\displaystyle \div_x \bigl(\A_{2}: \D(x,\frac{\partial
\w}{\partial t})+\A_{3}: \D(x,\w)+B_{4}\div_x \w +\\[1ex]
 \displaystyle
\int_0^t \bigl(\A_{4}(t-\tau ):\D(x,\w(\x,\tau )) +
 B_{5}(t-\tau )\div_x \w(\x,\tau )\bigr)d\tau \bigr),
\end{array} \right\}
\end{equation}
\begin{equation}\label{(2.36)}
\left. \begin{array}{r}\displaystyle \frac{1}{p_{*}}p+ m\div_x \w=\\[1ex]
\displaystyle -\int_0^t \bigl(C_{2}(t-\tau ):\D(x,\w(\x,\tau )) +
 a_{2}(t-\tau )\div_x \w(\x,\tau )\bigr)d\tau.
\end{array} \right\}
\end{equation}
\begin{equation}\label{(2.37)}
\left. \begin{array}{r}\displaystyle \frac{1}{\eta_{0}}\pi + (1- m)\div_x \w=\\[1ex]
\displaystyle-\int_0^t \bigl(C_{3}(t-\tau ):\D(x,\w(\x,\tau )) +
 a_{3}(t-\tau )\div_x \w(\x,\tau )\bigr)d\tau ,
\end{array} \right\}
\end{equation}
\begin{equation}\label{(2.38)}
q=p +\frac{\nu_0}{p_*}\frac{\partial p}{\partial t} .
\end{equation}

Here $\A_{2}$, $\A_{2}$ and  $\A_{4}$  are  fourth-rank tensors,
$B_{4}$, $B_{5}$,  $C_{2}$  and $C_{3}$ are matrices and $a_{2}$
and $a_{3}$ are scalars.  The exact expressions for these objects
are given below by formulas \eqref{(7.20)}--\eqref{(7.25)}.

The problem is supplemented by the homogeneous initial and
boundary conditions. (see in Eqs. \eqref{(2.12)} in Theorem
\ref{theorem2.2}).

 If the porous space is connected then  $\A_{2}$
  is strictly positively defined symmetric tensor.

If the porous space is disconnected, which is a case of isolated
pores, then   $\A_{2}=0$ and the system \eqref{(2.35)} degenerates
into nonlocal anisotropic Lam'e's system with strictly positively
defined and symmetric tensor  $\A_{3}$.
\end{theorem}

\addtocounter{section}{1} \setcounter{equation}{0}
\setcounter{theorem}{0} \setcounter{lemma}{0}
\setcounter{proposition}{0} \setcounter{corollary}{0}
\setcounter{definition}{0} \setcounter{assumption}{0}
\begin{center} \textbf{\S3. Preliminaries}
\end{center}

\textbf{3.1. Two-scale convergence.} Justification of Theorems
\ref{theorem2.1}--\ref{theorem2.4} relies on systematic use of the
method of two-scale convergence, which had been proposed by G.
Nguetseng \cite{NGU} and has been applied recently to a wide range
of homogenization problems (see, for example, the survey
\cite{LNW}).

\begin{definition} \label{TS}
A sequence $\{\varphi^\varepsilon\}\subset L^2(\Omega_{T})$ is
said to be \textit{two-scale convergent} to a limit $\varphi\in
L^2(\Omega_{T}\times Y)$ if and only if for any 1-periodic in $\y$
function $\sigma=\sigma(\x,t,\y)$ the limiting relation
\begin{equation}\label{(3.1)}
\lim_{\varepsilon\searrow 0} \int_{\Omega_{T}}
\varphi^\varepsilon(\x,t) \sigma\left(\x,t,\x /
\varepsilon\right)d\x dt = \int _{\Omega_{T}}\int_Y
\varphi(\x,t,\y)\sigma(\x,t,\y)d\y d\x dt
\end{equation}
holds.
\end{definition}

Existence and main properties of weakly convergent sequences are
established by the following fundamental theorem \cite{NGU,LNW}:
\begin{theorem} \label{theorem3.1}(\textbf{Nguetseng's theorem})

\textbf{1.} Any bounded in $L^2(Q)$ sequence contains a
subsequence, two-scale convergent to some limit
$\varphi\in L^2(\Omega_{T}\times Y)$.\\[1ex]
\textbf{2.} Let sequences $\{\varphi^\varepsilon\}$ and
$\{\varepsilon \nabla_x \varphi^\varepsilon\}$ be uniformly
bounded in $L^2(\Omega_{T})$. Then there exist a 1-periodic in
$\y$ function $\varphi=\varphi(\x,t,\y)$ and a subsequence
$\{\varphi^\varepsilon\}$ such that $\varphi,\nabla_y \varphi\in
L^2(\Omega_{T}\times Y)$, and $\varphi^\varepsilon$ and
$\varepsilon \nabla_x \varphi^\varepsilon$ two-scale converge to
$\varphi$ and $\nabla_y \varphi$,
respectively.\\[1ex]
\textbf{3.} Let sequences $\{\varphi^\varepsilon\}$ and
$\{\nabla_x \varphi^\varepsilon\}$ be bounded in $L^2(Q)$. Then
there exist functions $\varphi\in L^2(\Omega_{T})$ and $\psi \in
L^2(\Omega_{T}\times Y)$ and a subsequence from
$\{\varphi^\varepsilon\}$ such that $\psi$ is 1-periodic in $\y$,
$\nabla_y \psi\in L^2(\Omega_{T}\times Y)$, and
$\varphi^\varepsilon$ and $\nabla_x \varphi^\varepsilon$ two-scale
converge to $\varphi$ and $\nabla_x \varphi(\x,t)+\nabla_y
\psi(\x,t,\y)$, respectively.
\end{theorem}

\begin{corollary} \label{corollary3.1}
Let $\sigma\in L^2(Y)$ and
$\sigma^\varepsilon(\x):=\sigma(\x/\varepsilon)$. Assume that a
sequence $\{\varphi^\varepsilon\}\subset L^2(\Omega_{T})$
two-scale converges to $\varphi \in L^2(\Omega_{T}\times Y)$. Then
the sequence $\sigma^\varepsilon \varphi^\varepsilon$ two-scale
converges to $\sigma \varphi$.
\end{corollary}

\textbf{3.2. An extension lemma.} The typical difficulty in
homogenization problems while passing to a limit in Model
$B^\varepsilon$ as $\varepsilon \searrow 0$ arises because of the
fact that the bounds on the gradient of displacement $\nabla_x
\w^\varepsilon$ may be distinct in liquid and rigid phases. The
classical approach in overcoming this difficulty consists of
constructing of extension to the whole $\Omega$ of the
displacement field defined merely on $\Omega_s$. The following
lemma is valid due to the well-known results from \cite{ACE,JKO}.
We formulate it in appropriate for us form:

\begin{lemma} \label{Lemma3.1}
Suppose that assumptions of Sec. 1.2 on geometry of periodic
structure hold,   $ \psi^\varepsilon\in
W^1_2(\Omega^\varepsilon_s)$  and   $\psi^\varepsilon =0$ on
$S_{s}^{\varepsilon}=\partial\Omega ^\varepsilon_s \cap
\partial \Omega$ in the trace sense.  Then there exists a function
$ \sigma^\varepsilon \in
 W^1_2(\Omega)$ such that its restriction on the sub-domain
$\Omega^\varepsilon_s$ coincide with $\psi^\varepsilon$, i.e.,
\begin{equation} \label{3.2}
(1-\chi^\varepsilon(\x))( \sigma^\varepsilon(\x) -
\psi^\varepsilon (\x))=0,\quad \x\in\Omega,
\end{equation}
and, moreover, the estimate
\begin{equation} \label{3.3}
\|\sigma^\varepsilon\|_{2,\Omega}\leq C\|
\psi^\varepsilon\|_{2,\Omega ^{\varepsilon}_{s}}  , \quad
\|\nabla_x \sigma^\varepsilon\|_{2,\Omega} \leq  C \|\nabla_x
 \psi^\varepsilon\|_{2,\Omega ^{\varepsilon}_{s}}
\end{equation}
hold true, where the constant $C$ depends only on geometry $Y$ and
does not depend on $\varepsilon$.
\end{lemma}

\textbf{3.3. Friedrichs--Poincar\'{e}'s inequality in periodic
structure.} The following lemma was proved by L. Tartar in
\cite[Appendix]{S-P}. It specifies Friedrichs--Poincar\'{e}'s
inequality for $\varepsilon$-periodic structure.
\begin{lemma} \label{F-P}
Suppose that assumptions on the geometry of $\Omega^\varepsilon_f$
hold true. Then for any function $\varphi\in
\stackrel{\!\!\circ}{W^1_2}(\Omega^\varepsilon_f)$ the inequality
\begin{equation} \label{(F-P)}
\int_{\Omega^\varepsilon_f} |\varphi|^2 d\x \leq C \varepsilon^2
\int_{\Omega^\varepsilon_f} |\nabla_x \varphi|^2 d\x
\end{equation}
holds true with some constant $C$, independent of $\varepsilon$.
\end{lemma}

\textbf{3.4. Some notation.} Further we denote

 1) $$ \langle\Phi \rangle_{Y} =\int_Y \Phi  dy, \quad
 \langle\Phi \rangle_{Y_{f}} =\int_Y \chi \Phi  dy,
 \quad
 \langle\Phi \rangle_{Y_{s}} =\int_Y (1-\chi )\Phi  dy,$$
$$\langle\varphi  \rangle_{\Omega } =\int_{\Omega } \varphi  dx,
\quad
  \langle\varphi  \rangle_{\Omega_{T}} =\int_{\Omega_{T}} \varphi
  dxdt.$$
2) If $\textbf{a}$ and $\textbf{b}$ are two vectors then the
matrix $\textbf{a}\otimes \textbf{b}$ is defined by the formula
$$(\textbf{a}\otimes \textbf{b})\cdot
\textbf{c}=\textbf{a}(\textbf{b}\cdot \textbf{c})$$ for any vector
$\textbf{c}$.

3) If $B$ and $C$ are two matrices, then $B\otimes C$ is a
forth-rank tensor such that its convolution with any matrix $A$ is
defined by the formula
$$(B\otimes C):A=B (C:A)$$.
4) By $\I^{ij}$ we denote the $3\times 3$-matrix with just one
non-vanishing entry, which is equal to one and stands in the
$i$-th row and the $j$-th column.

5) We  also  introduce
$$
J^{ij}=\frac{1}{2}(\I^{ij}+\I^{ji})=\frac{1}{2} ({\mathbf e}_i
\otimes {\mathbf e}_j + {\mathbf e}_j \otimes {\mathbf e}_i),
$$
where $({\mathbf e}_1, {\mathbf e}_2, {\mathbf e}_3)$ are the
standard Cartesian basis vectors.

\addtocounter{section}{1} \setcounter{equation}{0}
\setcounter{theorem}{0} \setcounter{lemma}{0}
\setcounter{proposition}{0} \setcounter{corollary}{0}
\setcounter{definition}{0} \setcounter{assumption}{0}

\begin{center} \textbf{\S4. Proof of Theorem \ref{theorem2.1}}
\end{center}

\textbf{4.1.} Let  $\lambda_{0} <\infty$.

  If restriction $\tau_{0}>0$ holds, then
 \eqref{(2.2)}  follow from the lemma \ref{Lemma1.1}.

  Let $p_{*}<\infty$  and $\eta_{0}<\infty$.  Then pressures $p^{\varepsilon}$
  and  $\pi^{\varepsilon}$ are bounded from equations \eqref{(13)}
with the help of \eqref{(18)}. The pressure  $q^\varepsilon $ is
bounded from
 \begin{equation}\label{4.01}
q^\varepsilon =p^{\varepsilon} - \alpha
_{\nu}\chi^\varepsilon\div_x \frac{\partial \w}{\partial
t}^\varepsilon
\end{equation}
 with the help of  \eqref{(18)}.

If $p_{*}=\infty$, then  estimate \eqref{(2.4)} for the sum of
pressures $(q^\varepsilon+\pi^{\varepsilon})$ follows from the
basic integral identity \eqref{(15)}  and estimates
 \eqref{(18)} as an estimate for the corresponding functional,
  if we re-normalize the pressures  $(q^\varepsilon +\pi^{\varepsilon})$  such that
 $$\int _{\Omega}(q^\varepsilon(\x,t)+\pi^{\varepsilon}(\x,t))d\x=0. $$.
Indeed, the basic integral identity \eqref{(15)} and estimates
 \eqref{(18)} imply
$$|\int _{\Omega}( q^\varepsilon +\pi^{\varepsilon})\div_x
{\mathbf{\psi}} d\x|\leq C \|\nabla
{\mathbf{\psi}}\|_{2,\Omega}.$$
 Choosing now
${\mathbf{\psi}}$ such that  $(q^\varepsilon+\pi^\varepsilon )=
\div_x {\mathbf{\psi}}$ we get the desired estimate for the sum of
pressures  $(q^\varepsilon+\pi^\varepsilon )$. Such a choice is
always possible (see \cite{LAD}), if we put
$${\mathbf{\psi}}=\nabla \varphi + {\mathbf{\psi_{0}}}, \quad
\div_x {\mathbf{\psi_{0}}}=0, \quad \triangle
\varphi=q^\varepsilon+\pi^\varepsilon ,\quad \varphi |
_{\partial\Omega}=0, \quad (\nabla \varphi + {\mathbf{\psi_{0}}})|
_{\partial\Omega}=0.$$
 Note that the re-normalization of the pressures  $(q^\varepsilon+\pi^\varepsilon )$
  transforms continuity and state equations \eqref{(13)}
 for  pressures into
\begin{equation} \label{(4.1)}
\frac{1}{\alpha_p}p^\varepsilon +\chi^\varepsilon \div_x
\w^\varepsilon =-\frac{1}{m}\beta ^{\varepsilon}\chi^\varepsilon ,
 \end{equation}
 \begin{equation} \label{(4.1.1)}
\frac{1}{\alpha_\eta }\pi^\varepsilon +(1-\chi^\varepsilon )\div_x
\w^\varepsilon =(1-\chi^\varepsilon )\gamma ^{\varepsilon}_{s},
 \end{equation}
\begin{equation} \label{(4.1.2)}
q^\varepsilon =p^\varepsilon
+\frac{\alpha_{\nu}}{\alpha_p}\frac{\partial
p^\varepsilon}{\partial t}+\chi^\varepsilon \gamma
^{\varepsilon}_{f},
 \end{equation}
where
\begin{equation*}
\beta ^{\varepsilon}=\langle (1-\chi^\varepsilon )\div_x
\w^\varepsilon \rangle _{\Omega},\quad m\gamma
^{\varepsilon}_{f}=\langle q^\varepsilon \rangle _{\Omega},\quad
(1-m)\gamma ^{\varepsilon}_{s}=-\frac{1}{\alpha_\eta}\langle
q^\varepsilon \rangle _{\Omega}+\beta ^{\varepsilon}.
\end{equation*}

The case $\eta_{0}=\infty$ is considered in the same way. Note
that for all situations the basic integral identity \eqref{(15)}
permits to bound only the sum $(q^\varepsilon
+\pi^{\varepsilon})$. But thanks to the property that the product
of these two functions is equal to  zero, it is enough to get
bounds for each of these functions. The pressure $p^{\varepsilon}$
is bounded from the  state equation
 \eqref{(4.1.2)}, if we substitute the term  $(\alpha_{\nu} / \alpha_p)\partial
p^\varepsilon / \partial t$ from the state equation \eqref{(4.1)}
and use estimate \eqref{(18)}.

 Estimation of $\w^\varepsilon $  in the
case $\tau_0=0$ is not simple, and we outline it in more detail.

Let  $\mu_{1}>0$ and   $\tau_0=0$.  As usual, we obtain the basic
estimates if we multiply  equations for $\w^\varepsilon$
 by $\partial \w^\varepsilon / \partial t$ and then
integrate by parts all obtained terms. The only one term $\F\cdot
\partial \w^\varepsilon / \partial t$ heeds additional
consideration here. First of all, on the strength of Lemma
\ref{Lemma3.1}, we construct an extension $\uu^\varepsilon $ of
the function $\w^\varepsilon $ from $\Omega_s^\varepsilon$ into
$\Omega_f^\varepsilon $ such that $\uu^\varepsilon
=\w^\varepsilon$ in $\Omega_s^\varepsilon$, $\uu^\varepsilon \in
W_2^1(\Omega)$ and
$$\| \uu^\varepsilon\|_{2,\Omega} \leq C
\|\nabla_x \uu^\varepsilon\|_{2,\Omega} \leq
\frac{C}{\sqrt{\alpha_\lambda}}
 \|(1-\chi^\varepsilon)\sqrt{\alpha_\lambda}\nabla_x \w^\varepsilon\|_{2,\Omega }.$$

After that we estimate $\|\w^\varepsilon\|_{2,\Omega}$ with the
help of  Friedrichs--Poincar\'{e}'s inequality in periodic
structure (lemma \ref{F-P}) for the difference $(\uu^\varepsilon
-\w^\varepsilon)$:

$$\|\w^\varepsilon\|_{2,\Omega} \leq
\|\uu^\varepsilon\|_{2,\Omega} + \|\uu^\varepsilon
-\w^\varepsilon\|_{2,\Omega} \leq \|\uu^\varepsilon\|_{2,\Omega} +
C\varepsilon \|\chi^\varepsilon \nabla_x (\uu^\varepsilon
-\w^\varepsilon)\|_{2,\Omega} $$
$$\leq
\|\uu^\varepsilon\|_{2,\Omega}+C\varepsilon \|\nabla_x
\uu^\varepsilon\|_{2,\Omega}+C(\varepsilon \alpha _{\mu
}^{-\frac{1}{2}})\|\chi^\varepsilon \sqrt{\alpha_\mu} \nabla_x
\w^\varepsilon\|_{2,\Omega}$$
$$\leq \frac{C}{\sqrt{\alpha_\lambda}}
\|(1-\chi^\varepsilon)\sqrt{\alpha_\lambda}\nabla_x
\w^\varepsilon\|_{2,\Omega }+C(\varepsilon \alpha _{\mu
}^{-\frac{1}{2}})\|\chi^\varepsilon \sqrt{\alpha_\mu} \nabla_x
\w^\varepsilon\|_{2,\Omega}.$$
 Next we pass the derivative with respect to time from
$\partial \w^\varepsilon / \partial t$  to $\rho^{\varepsilon}\F$
and bound all positive terms (including the term $\alpha
_{\nu}\chi^\varepsilon\div_x \partial \w^\varepsilon /
\partial t $) in a usual way with
the help of H\"{o}lder and Grownwall's inequalities.

The rest of the proof is the same as for the case  $\tau_0>0$ if
we use a consequence of \eqref{(18)}:
$$\max\limits_{0<t<T}\alpha_\tau \| \frac{\partial ^2
\w^{\varepsilon}}{\partial t^2}(t)\|_{2,\Omega}\leq C.$$.\\
\noindent \textbf{4.2.} The proof of this part of the theorem is
obvious, because the re-normalization reduces this case to the
case of $\mu_{1}=1$ and  $\tau_0=0$, which has been already
considered.\\
\noindent \textbf{4.3.} Let  $\lambda_{0}=\infty$,
 $\mu_{1}<\infty$ and conditions \eqref{(2.9.1)} hold true. It is
 obvious that estimates  \eqref{(2.2)}  are still valid.

The desired estimates \eqref{(2.9.2)}  follow from the basic
equations for  $\alpha_\lambda \w^{\varepsilon}$ in the same way
as in the case of estimates \eqref{(2.2)}. The main difference
here is in the term
 $\rho ^{\varepsilon}\F\cdot\alpha_\lambda \partial \w^\varepsilon / \partial t$,
  which now  transforms to
$$\Upsilon \equiv \rho _{f}\F\cdot\alpha_\lambda
\frac{\partial \w^\varepsilon}{\partial t} +
 (\rho _{f}-\rho _{f})(1-\chi^\varepsilon )\F\cdot\alpha_\lambda
 \frac{\partial \w^\varepsilon}{\partial t}.$$
The integral of first term in $\Upsilon$ transforms as
  $$\rho _{f}\int_{0}^{t}\int_\Omega\nabla \Phi
  \cdot\alpha_\lambda \frac{\partial \w^\varepsilon}{\partial \tau}d\x d\tau=
 -\rho _{f}\int_{0}^{t}\int_\Omega \Phi \alpha_\lambda
 \div_x \frac{\partial \w^\varepsilon}{\partial \tau}d\x d\tau $$
 $$=-\rho _{f}\int_{0}^{t}\int_\Omega(\chi^\varepsilon
 \cdot\Phi \alpha_\lambda \div_x \frac{\partial \w^\varepsilon}{\partial \tau}
 +(1-\chi^\varepsilon)\cdot\Phi \alpha_\lambda \div_x
 \frac{\partial \w^\varepsilon}{\partial \tau})d\x d\tau$$
 $$=-\rho _{f}\int_\Omega(\chi^\varepsilon
 \cdot\Phi \alpha_\lambda \div_x \w^{\varepsilon}
 +(1-\chi^\varepsilon)\cdot\Phi \alpha_\lambda \div_x \uu^{\varepsilon})d\x +$$
 $$\rho _{f}\int_{0}^{t}\int_\Omega(\chi^\varepsilon
 \cdot\Phi _{\tau}\alpha_\lambda \div_x \w^{\varepsilon}
 +(1-\chi^\varepsilon)\cdot
 \Phi _{\tau}\alpha_\lambda \div_x \uu^{\varepsilon})d\x d\tau$$
and is bounded with the help of  terms
$$\int_\Omega (\chi^\varepsilon (\alpha _{p} \alpha_\lambda ^{-1})(\alpha_\lambda
\div_x \w^{\varepsilon})^{2}   + (1-\chi^\varepsilon)
|\alpha_\lambda \nabla_x \uu^\varepsilon |^{2})d\x,$$
 which appear in the basic identity after using the continuity equations.

 The integral of the second term in $\Upsilon$ is bounded with the
 help of the term
$$\langle (1-\chi^\varepsilon)|\alpha_\lambda \nabla_x
 \uu^\varepsilon |^{2} \rangle _{\Omega}$$ in the same way as before.

  Estimates \eqref{(2.4)} follow now from \eqref{(2.9.2)}.
Here as before the sum of  pressures $(q^\varepsilon +
\pi^\varepsilon)$ is bounded from the basic integral identity
\eqref{(15)} as a corresponding functional, and the pressure
$p^\varepsilon$ is bounded from the  equation \eqref{4.01}  due to
the bound \eqref{(18)} for  the divergence of the velocity of the
liquid component
 $\chi^\varepsilon \div_x (\partial \w^\varepsilon / \partial t)$.

If instead of conditions \eqref{(2.9.1)} one has condition
\eqref{(2.6)}, then bounds \eqref{(2.4)} for pressures
$p^\varepsilon$ and  $q^\varepsilon$  follow from equations
 \eqref{(13)} and  \eqref{4.01} and bounds \eqref{(18)}.
Note that in this case $\beta ^{\varepsilon}=0.$\qed

\addtocounter{section}{1} \setcounter{equation}{0}
\setcounter{theorem}{0} \setcounter{lemma}{0}
\setcounter{proposition}{0} \setcounter{corollary}{0}
\setcounter{definition}{0} \setcounter{assumption}{0}

\begin{center} \textbf{\S5. Proof of Theorem \ref{theorem2.2}}
\end{center}

\textbf{5.1. Weak and two-scale limits of sequences of
displacement and pressures.} On the strength of Theorem
\ref{theorem2.1}, the sequences $\{p^\varepsilon\}$,
$\{q^\varepsilon\}$, $\{\pi^\varepsilon\}$  and  $\{\w^\varepsilon
\}$   are uniformly in $\varepsilon$ bounded in $L^2(\Omega_{T})$.
Hence there exist a subsequence of small parameters
$\{\varepsilon>0\}$ and functions $p$, $q$, $\pi$ and  $\w$  such
that
\begin{equation}\label{(5.1)}
p^\varepsilon \rightarrow p,\quad q^\varepsilon \rightarrow q,
\quad \pi^\varepsilon \rightarrow \pi,  \quad  \w^\varepsilon
\rightarrow \w
\end{equation}
weakly in  $L^2(\Omega_T)$ as $\varepsilon\searrow 0$.

Due to Lemma \ref{Lemma3.1} there is a function $\uu^\varepsilon
\in L^\infty (0,T;W^1_2(\Omega))$ such that $\uu^\varepsilon
=\w^\varepsilon $ in $\Omega_{s}\times (0,T)$, and the family
$\{\uu^\varepsilon \}$ is uniformly in $\varepsilon$ bounded in
$L^\infty (0,T;W^1_2(\Omega))$. Therefore it is possible to
extract a subsequence of $\{\varepsilon>0\}$ such that
\begin{equation} \label{(5.2)}
\uu^\varepsilon \rightarrow \uu \mbox{ weakly in } L^2
(0,T;W^1_2(\Omega))
\end{equation}
as $\varepsilon \searrow 0$. Moreover,
\begin{equation} \label{(5.3)}
\chi^\varepsilon \alpha_\mu^\varepsilon \D_x(\w^\varepsilon)
\rightarrow 0.
\end{equation}

Relabelling if necessary, we assume that the sequences converge
themselves.

On the strength of Nguetseng's theorem, there exist 1-periodic in
$\y$ functions $P(\x,t,\y)$, $\Pi(\x,t,\y)$, $Q(\x,t,\y)$,
$\W(\x,t,\y)$  and $\UU(\x,t,\y)$ such that the sequences
$\{p^\varepsilon\}$, $\{\pi^\varepsilon\}$, $\{q^\varepsilon\}$,
$\{\w^\varepsilon \}$  and $\{\nabla_x \uu^\varepsilon \}$
two-scale converge to $P(\x,t,\y)$, $\Pi(\x,t,\y)$, $Q(\x,t,\y)$,
$\W(\x,t,\y)$  and $\nabla _{x}\uu +\nabla_{y}\UU(\x,t,\y)$, respectively.

Note that  the sequence  $\{\div_x \w^\varepsilon \}$ weakly
converges to $\div_x \w$ and $ \uu \in L^2
(0,T;\stackrel{\!\!\circ}{W^1_2}(\Omega)).$   Last assertion for
disconnected porous space follows from inclusion $\uu ^\varepsilon
\in L^2 (0,T;\stackrel{\!\!\circ}{W^1_2}(\Omega))$ and for the
connected porous space it follows from the
Friedrichs--Poincar\'{e}'s inequality for $\uu^\varepsilon$ in the
$\varepsilon$-layer  of the boundary $S$ and from convergence of
sequence $\{\uu^\varepsilon \}$  to $\uu$  strongly
 in $L^2(\Omega_{T})$ and weakly in $L^2 ((0,T);W^1_2(\Omega))$.\\

 \textbf{5.2. Micro- and macroscopic equations I.}
\begin{lemma} \label{lemma5.1}
For all $ \x \in \Omega$ and $\y\in Y$ weak  and two-scale limits
of the sequences $\{p^\varepsilon\}$, $\{\pi^\varepsilon\}$,
$\{q^\varepsilon\}$, $\{\w^\varepsilon\}$, and
$\{\uu^\varepsilon\}$ satisfy the relations
\begin{eqnarray} \label{(5.1)}
& P=\frac{1}{m}\chi p, \quad
Q=\frac{1}{m}\chi q;\\
\label{(5.2)} &  \frac{1}{\eta_{0}}\Pi+(1-\chi )
 (\div_x\uu + \div_y \UU)=\gamma _{s}(1-\chi );\\
 \label{(5.3)} & \div_y \W=0;\\
\label{(5.4)} &\W=\chi(\y) \W +(1-\chi)\uu;\\
\label{(5.5)} & q=p +\nu_0 p_{*}^{-1} \frac{\partial p}{\partial t}+\gamma _{f}m;\\
\label{(5.6)} & \frac{1}{p_{*}}p+\div_x \w = (1-m)\div_x \uu +
\langle \div_y\UU\rangle_{Y_{s}}-\beta ;\\
\label{(5.7)} & \frac{1}{\eta_{0}}\pi+(1-m)\div_x \uu + \langle
 \div_y\UU\rangle_{Y_{s}}=(1-m)\gamma
_{s},
\end{eqnarray}
where
$$\beta =\int_{\Omega}\langle \div_y\UU\rangle_{Y_{s}} dx,
\quad m\gamma _{f}=\langle q\rangle _{\Omega},\quad (1-m)\gamma
_{s}=-\frac{1}{\eta_{0}}\langle q\rangle _{\Omega}+\beta,$$  if
$p_{*}+\eta_{0}=\infty$  and $\beta=\gamma _{f}=\gamma _{s}=0$, if
$p_{*}+\eta_{0}<\infty.$
\end{lemma}

\begin{proof}
In order to prove Eq. \eqref{(5.1)}, into Eq. \eqref{(15)} insert
a test function ${\mathbf \psi}^\varepsilon =\varepsilon {\mathbf
\psi}\left(\x,t,\x / \varepsilon\right)$, where ${\mathbf
\psi}(\x,t,\y)$ is an arbitrary 1-periodic and finite on $Y_f$
function in $\y$. Passing to the limit as $\varepsilon \searrow
0$, we get
\begin{equation} \label{(5.8)}
\nabla_y Q(\x,t,\y)=0, \quad \y\in Y_{f}.
\end{equation}
The weak and two-scale limiting passage  in Eq. \eqref{(4.1.2)}
yield that Eq. \eqref{(5.5)} and the equation
\begin{equation} \label{(5.9)}
Q=P +\frac{\nu_0}{p_{*}} \frac{\partial P}{\partial t}+\gamma
_{f}\chi.
\end{equation}
hold. Taking into account Eq. \eqref{(5.8)} and \eqref{(5.9)}  we
get
\begin{equation*}
\nabla_y P(\x,t,\y)=0, \quad \y\in Y_{f}.
\end{equation*}
Next, fulfilling the two-scale limiting passage in  equalities
$$(1-\chi^{\varepsilon})p^{\varepsilon}=0, \quad (1-\chi^{\varepsilon})q^{\varepsilon} =0$$
we arrive at
$$(1-\chi )P=0, \quad (1-\chi )Q=0,$$
which along with Eqs. \eqref{(5.8)} and \eqref{(5.9)} justifies
Eq. \eqref{(5.1)}.

Eqs. \eqref{(5.2)}, \eqref{(5.3)}, \eqref{(5.6)}, and
\eqref{(5.7)} appear as the results of two-scale limiting passages
in Eqs. \eqref{(4.1)}-- \eqref{(4.1.2)} with the proper test
functions being involved. Thus, for example, Eq. \eqref{(5.6)}
arises if we represent Eq. \eqref{(4.1)} in the form
\begin{equation}\label{(5.10)}
p^\varepsilon + \div_x \w^\varepsilon =(1-\chi^\varepsilon )
\div_x \uu^\varepsilon -\frac{1}{m}\beta
^{\varepsilon}\chi^\varepsilon,
\end{equation}
multiply by an arbitrary function, independent of the ``fast''
variable $\x/\varepsilon$, and then pass to the limit as
$\varepsilon\searrow 0$. Eq. \eqref{(5.7)} is derived quite
similarly. In order to prove Eq. \eqref{(5.3)}, it is sufficient
to consider the two-scale limiting relations in Eq. \eqref{(5.10)}
 as $\varepsilon \searrow 0$ with the test functions $\varepsilon
\psi \left(\x / \varepsilon\right) h(\x,t)$, where $\psi$ and $h$
are arbitrary smooth test functions. In order to prove Eq.
\eqref{(5.4)} it is sufficient to consider the two-scale limiting
relations in
\begin{equation*}
(1-\chi ^{\varepsilon})(\w^{\varepsilon}-\uu^{\varepsilon})=0.
\end{equation*}
\end{proof}

\begin{lemma} \label{lemma5.2} For all $(\x,t) \in \Omega_{T}$
the relations
\begin{eqnarray} \label{(5.11)}
 &&\displaystyle \lambda_0 \triangle _{y}\UU
 = \nabla_y \Pi, \quad \y\in Y_s,\\
 \label{(5.12)}  && \displaystyle \bigl(\lambda_0
\D(y,\UU)-\Pi\cdot \I
 +\lambda_0 \D(x,\uu)
+ \frac{1}{m}q \cdot \I\bigr)\cdot {\mathbf n}=0, \quad \y\in
\gamma,
\end{eqnarray}
hold true. Here ${\mathbf n}$ is a unit normal to $\gamma$.
\end{lemma}

\begin{proof}
Substituting a test function of the form ${\mathbf
\psi}^\varepsilon =\varepsilon {\mathbf \psi}\left(\x,t,\x /
\varepsilon \right)$, where ${\mathbf \psi}(\x,t,\y)$ is an
arbitrary 1-periodic in $\y$ function vanishing on the boundary
$S$, into Eq.\eqref{(15)}, and passing to the limit as
$\varepsilon \searrow 0$, we arrive at the following microscopic
relation on the cell $Y$:
\begin{equation} \label{(5.14)}
\div_y \{\lambda_0(1-\chi ) (\D(y,\UU)+\D(x,\uu))- (\Pi
+\frac{1}{m}q \chi )\cdot \I \}=0,
\end{equation}
which is clearly equivalent to Eqs. \eqref{(5.11)} and
\eqref{(5.12)} in view of Eqs.\eqref{(5.1)}.
\end{proof}

\begin{lemma} \label{lemma5.3}
Let $\hat{\rho}=m \rho_{f} + (1-m)\rho_{s}$, $\V = \chi
\partial \w / \partial t$, and
 $\vv =\langle \V \rangle _{Y}$.
Then for all $0\leq \tau _{0}<\infty$ the quadruple of functions
$\{\uu , \vv, q, \pi \}$ satisfies in $\Omega_{T}$ the system of
macroscopic equations
\begin{eqnarray}\label{(5.13)}
    &&\tau
    _{0}\rho_{f}\frac{\partial \vv}{\partial t}+\tau
    _{0}\rho_{s}(1-m)\frac{\partial ^2\uu}{\partial t^2}-\hat{\rho}\F=\\
    &&\div_x \{\lambda _{0}((1-m)\D(x,\uu)+
    \langle \D(y,\UU)\rangle _{Y_{s}}
    )-(q+\pi )\cdot \I \}.\nonumber
\end{eqnarray}
\end{lemma}
\begin{proof}
Eqs. \eqref{(5.13)} arise as the limit of Eqs. \eqref{(15)} with
test functions being finite in $\Omega_T$ and independent of
$\varepsilon$.
\end{proof}

\textbf{5.3. Micro- and macroscopic equations II.}
\begin{lemma} \label{lemma5.4}
If $\mu_{1}=\infty$ then the weak and two-scale limits of
$\{\uu^\varepsilon\}$ and $\{\w^\varepsilon\}$ coincide.
\end{lemma}
\begin{proof}
In order to verify, it is sufficient to consider the difference
$(\uu^\varepsilon -\w^\varepsilon)$ and apply
Friedrichs--Poincar'{e}'s inequality, just like in the proof of
Theorem \ref{theorem2.1}.
\end{proof}

\begin{lemma} \label{lemma5.5}
Let $\mu_1 <\infty$. Then the weak and two-scale limits of
 $\{q^\varepsilon\}$ and $\{\chi^\varepsilon\w^\varepsilon\}$
 satisfy the microscopic relations
\begin{equation}\label{(5.16)}
\tau_{0}\rho_{f}\frac{\partial \V}{\partial t}-\rho_{f}\F=
\mu_{1}\triangle_y \V -\nabla_y R -\nabla_x q , \quad \y \in
Y_{f},
\end{equation}
\begin{equation}\label{(5.17)}
    \V=\frac{\partial \uu}{\partial t}, \quad \y \in \gamma
\end{equation}
in the case $\mu_{1}>0$, and relations
\begin{equation}\label{(5.18)}
    \tau_{0}\rho_{f}\frac{\partial \V}{\partial t}= -\nabla_y R
    -\nabla _{x} q
    +\rho_{f}\F, \quad \y \in Y_{f};
\end{equation}
    \begin{equation}\label{(5.19)}
    (\V-\frac{\partial \uu}{\partial t})\cdot{\mathbf n}=0, \quad \y \in \gamma
\end{equation}
in the case $\mu_{1}=0$.

In Eq. \eqref{(5.19)} ${\mathbf n}$ is the unit normal to
$\gamma$.
\end{lemma}

\begin{proof}
 Differential equations \eqref{(5.16)} and \eqref{(5.18)} follow
 as $\varepsilon\searrow 0$
 from integral equality \eqref{(15)} with the test function ${\mathbf
\psi}={\mathbf \varphi}(x\varepsilon^{-1})\cdot h({\mathbf x},t)$,
where ${\mathbf \varphi}$ is solenoidal and finite in $Y_{f}$.

Boundary conditions \eqref{(5.17)} are the consequences of the
two-scale convergence of $\{\alpha_{\mu}^{\frac{1}{2}}\nabla_x
\w^{\varepsilon}\}$ to the function
$\mu_{1}^{\frac{1}{2}}\nabla_y\W(\x,t,\y)$. On the strength of
this convergence, the function $\nabla_y \W (\x,t,\y)$ is
$L^2$-integrable in $Y$. The boundary conditions \eqref{(5.19)}
follow from Eq. \eqref{(5.3)}.
\end{proof}

\begin{lemma} \label{lemma5.6}
If the porous space is disconnected, which is the case of isolated
pores, then the weak and two-scale limits of sequences
$\{\uu^\varepsilon\}$ and $\{\w^\varepsilon\}$ coincide.
\end{lemma}
\begin{proof}
Indeed, in the case $0\leq \mu_{1}<\infty$ the systems of
equations \eqref{(5.3)}, \eqref{(5.16)}, and \eqref{(5.17)}, or
\eqref{(5.3)}, \eqref{(5.18)}, and \eqref{(5.19)} have the unique
solution $\V=\partial \uu / \partial t$.
\end{proof}

\textbf{5.4. Homogenized equations I.}

\begin{lemma} \label{lemma5.7}
If $\mu_1 =\infty$ or the porous space is disconnected then
$\w=\uu$ and  the weak limits $\uu$, $p$, $q$, and $\pi$ satisfy
in $\Omega_{T}$ the initial-boundary value problem
 \begin{equation}\label{(5.20)}
\tau _{0}\hat{\rho}\frac{\partial ^2\uu}{\partial t^2}=
    \div_x \{\lambda _{0}\A^{s}_{0}:\D(x,\uu) +  B^{s}_{0}\div_x \uu
 +B^{s}_{1}q - (q+\pi )\cdot \I \}+\hat{\rho}\F,
\end{equation}
\begin{equation}\label{(5.21)}
\frac{1}{\eta_{0}}\pi+C^{s}_{0}:\D(x,\uu)+ a^{s}_{0}\div_x \uu
+a^{s}_{1}q+a^{s}_{2}\langle q\rangle _{\Omega}=0,
\end{equation}
\begin{equation}\label{(5.22)}
 q=p +\nu_0 p_{*}^{-1} \frac{\partial p}{\partial t}+ \gamma _{f}m,\quad
 \frac{1}{p_{*}}p + \frac{1}{\eta_{0}}\pi + \div_x \uu=(1-m)\gamma _{s}-\beta,
\end{equation}
where the symmetric strictly  positively defined constant
fourth-rank tensor $\A^{s}_{0}$, matrices  $C^{s}_{0}, B^{s}_{0}$
and  $B^{s}_{1}$ and constants $a^{s}_{0}$, $a^{s}_{1}$ and
$a^{s}_{2}$  are defined below by formulas \eqref{(5.26)},
\eqref{(5.27.1)} - \eqref{(5.27.3)}.

Differential equations \eqref{(5.20)} are endowed with homogeneous
initial and boundary conditions
 \begin{equation}\label{(5.22.1)}
 \tau _{0}\uu(\x,0)= \tau _{0}\frac{\partial \uu}{\partial t}(\x,0)=0,\quad \x\in \Omega,
 \quad \uu(\x,t)=0, \quad \x\in S, \quad t>0.
\end{equation}
\end{lemma}
\begin{remark}\label{remark1}
In what follows we can neglect terms $a^{s}_{2}\langle q\rangle
_{\Omega}$, $\gamma _{f}$ and  $\gamma _{s}$ because all pressures
are defined up to arbitrary functions of time.
\end{remark}
\begin{proof}
In the first place let us notice that $\uu =\w$ due to Lemmas
\ref{lemma5.4} and \ref{lemma5.6}.

The homogenized equations \eqref{(5.20)} follow from the
macroscopic equations \eqref{(5.13)}, after we insert in them the
expression
$$\langle \D(y,\UU)\rangle _{Y_{s}}=\A^{s}_{1}:\D(x,\uu) + B^{s}_{0}\div_x \uu
 +B^{s}_{1}q.$$
In turn, this expression follows by virtue of solutions of Eqs.
 \eqref{(5.2)} and \eqref{(5.14)} on the pattern cell $Y_{s}$.
 Indeed, setting
 $$ \UU=\sum_{i,j=1}^{3}\UU^{ij}(\y)D_{ij}+
 \UU_{0}(\y)\div_x \uu + \frac{1}{m}(\UU_{1}(\y)(q-\langle q\rangle _{\Omega})
 +\UU_{2}(\y)\langle q\rangle _{\Omega}), $$
 $$\Pi =\lambda _{0}\sum_{i,j=1}^{3}\Pi^{ij}(\y)D_{ij}
 +\Pi_{0}(\y)\div_x \uu + \frac{1}{m}(\Pi_{1}(\y)(q-\langle q\rangle
 _{\Omega})+\Pi_{2}(\y)\langle q\rangle _{\Omega}),$$

where
 $$D_{ij}(\x,t)=\frac{1}{2}(\frac{\partial u_{i}}{\partial x_{j}}(\x,t)+
 \frac{\partial u_{j}}{\partial x_{i}}(\x,t)),$$
we arrive at the following periodic-boundary value problems in
$Y_{s}$:
\begin{equation}\label{(5.25.1)}
\left. \begin{array}{lll}  \displaystyle \div_y \{(1-\chi )
(\D(y,\UU^{ij})+J^{ij}) - \Pi ^{ij}\cdot \I \}=0,\\[1ex]
\frac{\lambda _{0}}{\eta_{0}}\Pi ^{ij} +(1-\chi ) \div_y \UU^{ij}
=0;
\end{array} \right\}
\end{equation}
\begin{equation}\label{(5.25.2)}
\left. \begin{array}{lll}  \displaystyle \div_y
\{\lambda_{0}(1-\chi ) \D(y,\UU_{0}) - \Pi_{0}\cdot \I \}=0,\\[1ex]
\frac{1}{\eta_{0}}\Pi _{0} + (1-\chi )(\div_y \UU_{0}+1) =0;
\end{array} \right\}
\end{equation}
\begin{equation}\label{5.25.3}
\left. \begin{array}{lll}  \displaystyle \div_y
\{\lambda_{0}(1-\chi ) \D(y,\UU_{1}) - (\Pi_{1}+\chi )\cdot \I
\}=0,\\[1ex]
\frac{1}{\eta_{0}}\Pi _{1} +(1-\chi )\div_y \UU_{1}) =0.
\end{array} \right\}
\end{equation}
\begin{equation}\label{(5.25.4)}
\left. \begin{array}{r} \displaystyle \div_y \{\lambda_{0}(1-\chi
)\D(y,\UU_{2}) - (\Pi_{2}+\chi )\cdot \I \}=0,  \\[1ex]
\displaystyle \frac{1}{\eta_{0}}\Pi _{2} + (1-\chi )\div_y \UU_{2}
-\frac{(1-\chi )}{(1-m)}(\langle
\div_y\UU_{2}\rangle_{Y_{s}}+\frac{m}{\eta_{0}})=0.
\end{array} \right\}
\end{equation}
Note that for  $p_{*}+\eta_{0}=\infty$
  $$\beta=\langle\langle \div_y \UU \rangle_{\Omega}\rangle_{Y_{s}}=
  \sum_{i,j=1}^{3}\langle \div_y\UU^{ij}\rangle_{Y_{s}}
  \langle D_{ij}\rangle_{\Omega} +\langle \div_y\UU_{0}\rangle_{Y_{s}}
  \langle \div_x \uu\rangle_{\Omega} + $$
  $$\frac{1}{m}\langle \div_y\UU_{1}\rangle_{Y_{s}}
  \langle (q-\langle q\rangle_{\Omega})\rangle_{\Omega}+
  \frac{1}{m}\langle \div_y\UU_{2}\rangle_{Y_{s}}
  \langle q\rangle_{\Omega}=\frac{1}{m}\langle \div_y\UU_{2}\rangle_{Y_{s}}
  \langle q\rangle_{\Omega}$$
 due to homogeneous boundary conditions for   $\uu(\x,t)$.

On the strength of the assumptions on the geometry of the pattern
``liquid'' cell $Y_{s}$, problems  \eqref{(5.25.1)}--
\eqref{(5.25.3)} have  unique solution, up to an arbitrary
constant vector. In order to discard the arbitrary constant
vectors we demand
$$
\langle\UU^{ij}\rangle _{Y_{s}}=\langle\UU_{0}\rangle_{Y_{s}}
=\langle\UU_{1}\rangle_{Y_{s}} =\langle\UU_{2}\rangle_{Y_{s}}=0.
$$
Thus
 \begin{equation}\label{(5.26)}
 \A^{s}_{0}=\sum_{i,j=1}^{3}J^{ij}\otimes J^{ij} + \A^{s}_{1}, \quad
 \A^{s}_{1}=\sum_{i,j=1}^{3}\langle (1-\chi) D(y,\UU^{ij})\rangle _{Y}\otimes
    J^{ij}.
\end{equation}
 Symmetry of the tensor $\A^{s}_{0}$ follows from symmetry of
 the tensor $\A^{s}_{1}$. And symmetry of the latter one follows
 from the equality
 \begin{eqnarray}\label{(5.27)}
     &&\langle\D(y,\UU^{ij})\rangle _{Y_{s}} : J^{kl}
 =\\
 &&-\langle \D(y,\UU^{ij}) : \D(y,\UU^{kl})\rangle_{Y_{s}} +
 \frac{\lambda _{0}}{\eta_{0}}\Pi ^{ij} \Pi ^{kl},\nonumber
\end{eqnarray}
which appears by means of multiplication of Eq. \eqref{(5.25.1)}
for $\UU^{ij}$ by $\UU^{kl}$ and by integration by parts.

This equality also implies positive definiteness of the tensor
$\A^{s}_{0}$. Indeed, let $\zeta$ be an arbitrary symmetric
matrix. Setting
 $$\Z=\sum_{i,j=1}^{3}\UU^{ij}\zeta_{ij},\quad \tilde{\Pi}=
 \sum_{i,j=1}^{3}\Pi^{ij}\zeta_{ij}$$
and taking into account Eq. \eqref{(5.27)} we get
  \begin{equation*}
     \langle\D(y,\Z)\rangle _{Y_{s}}:\zeta
 =-\langle\D(y,\Z): \D(y,\Z)\rangle_{Y_{s}} -
 \frac{\lambda _{0}}{\eta_{0}}\tilde{\Pi}^{2},
\end{equation*}
 This equality and the definition of the tensor $A_{0}^s$ give
\begin{equation*}
 (A_{0}^s:\zeta ):\zeta =
 \langle(\D(y,\Z)+\zeta ): (\D(y,\Z)+\zeta )\rangle_{Y_{s}}
 + \frac{\lambda _{0}}{\eta_{0}}\tilde{\Pi}^{2}.
\end{equation*}
Now the strict positive definiteness of the tensor $\A_{0}$
follows from the equality immediately above and the geometry of
the elementary cell $Y_{s}$. Namely, let for some function $\zeta$
such that $\zeta :\zeta =1$ be $(\A_{0}:\zeta ):\zeta =0$. But
then we have $(\D(y,\Z)+\zeta )=0$, which is possible iff $\Z$ is
a linear function in $\y$. On the other hand, all linear periodic
functions on $Y_{s}$ are constant. Finally, the normalization
condition $\langle\UU^{ij}\rangle =0$ yields that $\Z=0$. However,
this is impossible because the functions $\UU^{ij}$ are linearly
independent.

 Finally, Eqs. \eqref{(5.21)} and \eqref{(5.22)} for the pressures follow from Eqs.
 \eqref{(5.5)}-- \eqref{(5.7)} and
\begin{equation}\label{(5.27.1)}
  B^{s}_{0}=\langle\D(y,\UU_{0})\rangle _{Y_{s}}, \quad
   B^{s}_{1}=\frac{1}{m}\langle\D(y,\UU_{1})\rangle _{Y_{s}},
\end{equation}
\begin{equation}\label{(5.27.2)}
 C^{s}_{0}=\sum_{i,j=1}^{3}\langle \div_y\UU^{ij}\rangle _{Y_{s}}J^{ij}, \quad
  a^{s}_{0}=1-m + \langle\div_y\UU_{0}\rangle _{Y_{s}},
\end{equation}
\begin{equation}\label{(5.27.3)}
\ a^{s}_{1}= \frac{1}{m}\langle\div_y\UU_{1}\rangle _{Y_{s}},\quad
   a^{s}_{2}= \frac{1}{m}(\frac{m}{\eta _{0}}-\langle\div_y\UU_{1}\rangle _{Y_{s}}).
\end{equation}
\end{proof}

\textbf{5.5. Homogenized equations II.}

Let $\mu_{1}<\infty$. In the same manner as above, we verify that
the limit $\uu$ of the sequence $\{\uu^\varepsilon\}$ satisfies
the initial-boundary value problem likes \eqref{(5.20)}--
\eqref{(5.22)}. The main difference here that, in general, the
weak limit $\w$ of the sequence $\{\W^\varepsilon\}$ differs from
$\uu$. More precisely, the following statement is true.
\begin{lemma} \label{lemma5.8}
If $\mu_{1}<\infty$ then the weak limits $\uu$, $\w^{f}$, $p$,
$q$, and $\pi$ of the sequences
 $\{\uu^\varepsilon\}$, $\{\chi^{\varepsilon}\W^\varepsilon\}$,
  $\{p^\varepsilon\}$,  $\{q^\varepsilon\}$, and $\{\pi^\varepsilon\}$
satisfy the initial-boundary value problem in $\Omega_T$,
consisting of the balance of momentum equation
\begin{eqnarray}\label{(5.28)}
&&\tau _{0}(\rho_{f}\frac{\partial \vv}{\partial t}+\rho_{s}
(1-m)\frac{\partial ^2\uu}{\partial t^2}) + \nabla (q+\pi )-\hat{\rho}\F= \\
&&\div_x \{\lambda _{0}A^{s}_{0}:\D(x,\uu) +
 B^{s}_{0}\div_x \uu
 +B^{s}_{1}q \},\nonumber
\end{eqnarray}
where   $\vv=\partial \w^{f} / \partial t$  and $\A^{s}_{0}$,
$B^{s}_{0}$  and $B^{s}_{1}$ are the same as in \eqref{(2.10)},
the continuity equations \eqref{(5.9)}, the equations
\begin{equation} \label{(5.28.1)}
     p +\nu_0 p_{*}^{-1} \frac{\partial p}{\partial t}= q,\quad
     \frac{1}{p_{*}}p+ \frac{1}{\eta_{0}}\pi+\div_x \w^{f}  =  (m-1)\div_x \uu,
   \end{equation}
 and Darcy's law in the form
\begin{equation}\label{(5.29)}
\vv=m\frac{\partial \uu}{\partial t}+\int_{0}^{t}
B_{1}(\mu_1,t-\tau)\cdot (-\nabla_x
q+\rho_{f}\F-\tau_{0}\rho_{f}\frac{\partial ^2 \uu}{\partial \tau
^2})(\x,\tau )d\tau
\end{equation}
in the case of $\tau_{0}>0$ and $\mu_{1}>0$, Darcy's law in the
form
\begin{equation}\label{(5.30)}
\vv=m\frac{\partial \uu}{\partial t}+B_{2}(\mu_1)\cdot(-\nabla_x
q+\rho_{f}\F)
\end{equation}
in the case of $\tau_{0}=0$ and $\mu_{1}>0$, and, finally, Darcy's
law in the form
\begin{equation}\label{(5.31)}
\vv=B_{3}\cdot \frac{\partial \uu}{\partial t}+\frac{1}{\tau
_{0}\rho_{f}}(m\I-B_{3})\cdot\int_{0}^{t}(-\nabla_x q(\x,\tau
)+\rho_{f}\F(\x,\tau ))d\tau
\end{equation}
in the case of $\tau_{0}>0$ and $\mu_{1}=0$. The problem is
supplemented by boundary and initial conditions \eqref{(5.22)} for
the displacement $\uu$ of the rigid component and by the boundary
condition
\begin{equation}\label{(5.32)}
\vv(\x,t)\cdot \n(\x)=0,
     \quad (\x,t) \in S=\partial \Omega , \quad t>0,
\end{equation}
for the velocity $\vv$ of the liquid component. In Eqs.
\eqref{(5.29)}--\eqref{(5.32)} $\n(\x)$ is the unit normal vector
to $S$ at a point $\x \in S$, and matrices $B_{1}(\mu_1,t)$,
$B_{2}(\mu_1)$, and $B_{3}$ are given below by Eqs.
\eqref{(5.33)}--\eqref{(5.38)}.
\end{lemma}
\begin{proof}
The derivation of Eq. \eqref{(5.32)} is standard \cite{S-P}. The
homogenized equations of balance of momentum and balance of mass
derive exactly as \eqref{(5.20)}--\eqref{(5.21)}. For example, to
get Eq. \eqref{(5.28.1)} we just expressed $\div_x \w $ in a sum
of Eqs. \eqref{(5.6)} and \eqref{(5.7)} using Eq. \eqref{(5.4)}
after homogenization: $\w=\w^{f}+(1-m)\uu.$ Therefore we omit the
relevant proofs now and focus only on derivation of homogenized
equations for the velocity $\vv$ in the form of Darcy's laws.

a) If  $\mu_{1}>0$ and  $\tau_{0}>0$, then the solution of the
system of microscopic equations \eqref{(5.3)}, \eqref{(5.16)}, and
\eqref{(5.17)}, provided with the homogeneous initial data, is
given by formula
\begin{equation*}
\V=\frac{\partial \uu}{\partial t}+\int_{0}^{t}
\textbf{B}^{f}_{1}(\y,t-\tau)\cdot (-\nabla_x
q+\rho_{f}\F-\tau_{0}\rho_{f}\frac{\partial ^2 \uu}{\partial \tau
^2})(\x,\tau )d\tau ,
\end{equation*}
 in which
\begin{equation*}
\textbf{B}^{f}_{1}(\y,t)= \sum_{i=1}^{3}\V^{i}(\y,t)\otimes
\e_{i},
\end{equation*}
and the functions $\V^{i}(\y,t)$ are defined by virtue of the
periodic initial-boundary value problem
\begin{equation}\label{(5.34)}
\left. \begin{array}{lll}  \displaystyle \tau
_{0}\rho_{f}\frac{\partial \V^{i}}{\partial t}-\mu_{1}\triangle
\V^{i} +\nabla Q^{i} =0,
  \quad \mbox{div}_y \V^{i} =0, \quad \y \in Y_{f},  t>0,\\[1ex]
\V^{i}=0, \quad \y \in \gamma ,  t>0;\quad
 \tau _{0}\rho_{f}\V^{i}(y,0)=\e_{i}, \quad \y \in Y_{f}.
 \end{array} \right\}
\end{equation}
 In Eq. \eqref{(5.34)} $\e_{i}$ is the standard Cartesian
 basis vector.

Therefore
\begin{equation}\label{(5.33)}
B_{1}(\mu_{1},t)= \langle \textbf{B}^{f}_{1}(\y,t)\rangle
_{Y_{s}},
\end{equation}
 b) If $\tau_{0}=0$ and $\mu_{1}>0$ then the solution of
the stationary microscopic equations \eqref{(5.3)},
\eqref{(5.16)}, and \eqref{(5.17)} is given by formula
\begin{equation*}
\V=\frac{\partial \uu}{\partial
t}+\textbf{B}^{f}_{2}(\y)\cdot(-\nabla q+\rho_{f}\F),
\end{equation*}
in which
\begin{equation*}
\textbf{B}^{f}_{2}(\y)= \sum_{i=1}^{3}\UU^{i}(\y)\otimes \e_{i} ,
\end{equation*}
and the functions $\UU^{i}(\y)$ are defined from the periodic
boundary value problem
\begin{equation}\label{(5.36)}
\left. \begin{array}{lll}  \displaystyle -\mu_{1}\triangle \UU^{i}
+\nabla R^{i} =\e_{i}, \quad \mbox{div}_y \UU^{i} =0, \quad \y \in
Y_{f},\\[1ex]
\UU^{i}=0, \quad \y \in \gamma .
\end{array} \right\}
\end{equation}
Thus
\begin{equation}\label{(5.35)}
B_{2}(\mu_{1})= \langle \textbf{B}^{f}_{2}((\y)\rangle _{Y_{s}}.
\end{equation}
 Matrices $B_{1}(\mu_1,t)$ and $B_{2}(\mu_1)$ are symmetric
and positively defined \cite[Chap. 8]{S-P}.

c) If $\tau_{0}>0$ and $\mu_{1}=0$ then in the process of solving
the system \eqref{(5.3)}, \eqref{(5.18)}, and \eqref{(5.19)} we
firstly find the pressure $R(\x,t,\y)$ by virtue of solving the
Neumann problem for Laplace's equation in
 $Y_{f}$. If
 $$
 \h(\x,t)=-\tau_{0}\rho_{f}\frac{\partial ^2\uu}{\partial t^2}(\x,t)
 -\nabla_x q(\x,t)+\rho_{f}\F(\x,t),
 $$
 then
 $$R(\x,t,\y)=\sum_{i=1}^{3} R_{i}(\y) \e_{i}\otimes \h(\x,t),$$
 where $R^{i}(\y)$ is the solution of the problem
 \begin{equation}\label{(5.37)}
\triangle_y R_{i}=0,\quad \y \in Y_{f}; \quad \nabla_y R_{i}\cdot
\n =\n\cdot \e_{i}, \quad \y \in \gamma .
\end{equation}
Formula \eqref{(5.31)} appears as the result of homogenization of
Eqs. \eqref{(5.18)} and
 \begin{equation}\label{(5.38)}
B_{3}=\sum_{i=1}^{3}\langle \nabla R_{i}(\y)\rangle
_{Y_{s}}\otimes \e_{i},
\end{equation}
where the matrix $(m\I - B_3)$ is symmetric and positively
definite \cite[Chap. 8]{S-P}.
\end{proof}

\addtocounter{section}{1} \setcounter{equation}{0}
\setcounter{theorem}{0} \setcounter{lemma}{0}
\setcounter{proposition}{0} \setcounter{corollary}{0}
\setcounter{definition}{0} \setcounter{assumption}{0}
\begin{center} \textbf{\S6. Proof of Theorem \ref{theorem2.3}}
\end{center}

\textbf{6.1. Weak and two-scale limits of sequences of
displacement and pressures.}

$\textbf{I.}$ Let $\mu_{1}<\infty$  and one of the conditions
  \eqref{(2.6)} or  \eqref{(2.9.1)}  holds true.
 Then on the strength of Theorems
\ref{theorem2.1} and \ref{theorem3.1} we conclude that sequences
 $\{\chi^{\varepsilon}\w^\varepsilon\}$,
  $\{p^\varepsilon\}$ and $\{q^\varepsilon\}$
two-scale converge to  $\chi (\y)\W(\x,t,\y)$,
  $P(\x,t,\y)$  and  $Q(\x,t,\y)$  and weakly converge in
$L^2(\Omega_{T})$ to  $\w^{f}$, $p$  and  $q$  respectively,  and
a sequence  $\{\uu^\varepsilon (\x,t)\}$, where  $\uu^\varepsilon
 (\x,t)$  is an extension of  $\w^\varepsilon (\x,t)$  from the
 domain $\Omega_{s}^{\varepsilon}$  into domain $\Omega$,
 strongly  converges in $L^2(\Omega_{T})$ and weakly in
$L^2((0,T);W^1_2(\Omega))$  to zero.

$\textbf{II.}$ If $\mu_{1}<\infty$ and conditions \eqref{(2.9.1)}
hold true, then due to estimates  \eqref{(2.4)} and
\eqref{(2.9.2)} the sequence $\{\alpha_{\lambda}\uu^\varepsilon\}$
 converges strongly in   $L^2(\Omega_{T})$ and weakly in
$L^2((0,T);W^1_2(\Omega))$ to a function $\uu$, and the sequence
$\{\pi^{\varepsilon}\}$   converges weakly in  $L^{2}(\Omega_{T})$
to a function $\pi$.

 $\textbf{III.}$ If  $\mu_{1}=\infty$,  $p_{1}^{-1}, \eta
_{1}^{-1}<\infty $  and  $0< \lambda_{1}<\infty$,
     then on the strength of part 2) of Theorem  \ref{theorem2.1}
     the sequences
   $\{\alpha_{\mu}\varepsilon^{-2}\chi^\varepsilon \w^\varepsilon\}$,
  $\{p^\varepsilon\}$, $\{\pi^\varepsilon\}$ and  $\{q^\varepsilon\}$
two-scale converge to functions  $\chi (\y)\W(\x,t,\y)$,
  $P(\x,t,\y)$, $\Pi(\x,t,\y)$ and  $Q(\x,t,\y)$ and weakly in
$L^2(\Omega_{T})$ to functions  $\w^{f}$, $p$, $\pi$ and  $q$
respectively, and the sequence
 $\{\alpha_{\mu}\varepsilon^{-2}\uu^\varepsilon \}$
strongly converge in $L^2(\Omega_{T})$ and weakly in
$L^2((0,T);W^1_2(\Omega))$ to the function $\uu$.

As before in previous section \S5, we conclude  that $ \uu \in L^2
(0,T;\stackrel{\!\!\circ}{W^1_2}(\Omega)).$

\textbf{6.2.Homogenized equations.}

  $\textbf{I.}$
If $\mu_{1}<\infty$ and one of the conditions  \eqref{(2.6)} or
\eqref{(2.9.1)}  holds true, then, as in the proof of Theorem
\ref{theorem2.2}, we construct a closed system of equations for
the velocity $\vv=\partial \w^{f} / \partial t$  in the liquid
component  and for
 the pressures $p$ and $q$, consisting of the modifications of Darcy's law
\eqref{(5.29)}--\eqref{(5.31)} and boundary condition
\eqref{(5.32)}, in which we have $\uu(\x,t)=0$, and of the
equations
\begin{equation} \label{(6.1)}
     p +\nu_0 p_{*}^{-1} \frac{\partial p}{\partial t}= q,\quad
     \frac{1}{p_{*}}\frac{\partial p}{\partial t}+\div_x \vv  = 0,
    \quad \x \in \Omega , \quad t>0.
\end{equation}

We entitle the above described systems as  Problem
 $ F_1$, $ F_2$, or $F_3$ depending on the forms of the
 matrices $B_1$, $B_2$, or $B_3$, having places in Darcy's laws.

$\textbf{II.}$  Let  $\mu_{1}<\infty$ and condition
\eqref{(2.9.1)} holds true. We observe that the limiting
displacements in the rigid skeleton are equal to zero. In order to
find a more accurate asymptotic of the solution of the original
model, we use again the re-normalization. Namely, let
 $$\w^{\varepsilon}\rightarrow\alpha_{\lambda}\w^{\varepsilon}.$$

Then new displacements satisfy the same problem as displacements
before re-normalization, but with new parameters
  $$\alpha_{\eta}\rightarrow
\alpha_{\eta}\alpha_{\lambda}^{-1}, \quad
 \alpha_{\lambda}\rightarrow 1, \quad \alpha_{\tau}\rightarrow
 \alpha_{\tau}\alpha_{\lambda}^{-1}.$$

 Thus we arrive at the assumptions of Theorem \ref{theorem2.2}.
Namely, the limiting functions $\uu(\x,t)$, $\pi(\x,t)$,
$\Pi(\x,t,\y)$, and $\UU(\x,t,\y)$ satisfy the system of micro-
and macroscopic equations \eqref{(5.2)}, \eqref{(5.14)},
\eqref{(5.7)} and
  \eqref{(5.13)}, in which the pressure $q$ is given by
virtue of one of Problems
 $ F_1$-- $ F_3$.
The only difference from already considered case is in micro- and
macroscopic continuity equations, because this equation depends on
the value $\eta_{2}$. These micro- and macroscopic continuity
equations coincide with Eq. \eqref{(5.2)} and  Eq. \eqref{(5.7)}
if we put there  $\eta_{0}=\eta_{2}$.

Hence for $\uu(\x,t)$ and $\pi(\x,t)$  there hold true the
homogenized momentum equation in the form
    \begin{equation}\label{(6.3)}
0=\div_x \{\A^{s}_{0}:\D(x,\uu) + B^{s}_{0}\div_x \uu
 +B^{s}_{1}q - (q+\pi )\cdot \I \}+\hat{\rho}\F, \quad \x \in \Omega ,
\end{equation}
continuity equation \eqref{(5.21)},  in which we
  have  $\eta_{0}=\eta_{2}$,
and the boundary condition  \eqref{(2.12)}.

 The tensor $\A^{s}_{0}$, the  matrices  $C^{s}_{0}, B^{s}_{0}$ and  $B^{s}_{1}$
and the constants $a^{s}_{0}$ и $a^{s}_{1}$  are defined  from
Eqs. \eqref{(5.26)}, \eqref{(5.27.1)}--\eqref{(5.27.2)}, in which
 we  have $\eta_{0}=\eta_{2}$ and $\lambda_{0}=1$.

 $\textbf{III.}$
If $\mu_{1}=\infty$, $p_{1}^{-1}, \eta _{1}^{-1}<\infty $  and
$0<\lambda_{1}<\infty$ then re-normalizing by
   $$\w^{\varepsilon} \rightarrow
   \alpha_{\mu}\varepsilon^{-2}\w^{\varepsilon}$$
  we arrive at the assumptions of Theorem \ref{theorem2.2}, when
   $\mu_{1}=1$, $\tau_{0}=0$ и $\lambda_{0}=\lambda_{1}$.
Namely, functions $\w^{f}$, $p$, $\pi$  and $\uu$ satisfy the
following initial-boundary value problem in $\Omega_T$:
\begin{equation}\label{(6.5)}
\left. \begin{array}{r} \displaystyle \div_x \{\lambda
_{1}\A^{s}_{0}:\D(x,\uu) + B^{s}_{0}\div_x \uu
 +B^{s}_{1}p - (p+\pi )\cdot \I \}+\hat{\rho}\F=0,\\[1ex]
\displaystyle \frac{\partial \w^{f}}{\partial t}=\frac{\partial
\uu}{\partial t}+B_{2}(1)\cdot (-\nabla p + \rho_{f}\F),\\[1ex]
\displaystyle \frac{1}{p_{1}} p+ \frac{1}{\eta _{1}}\pi + \div_x
\w^{f} =(m-1)\div_x \uu,\\[1ex]
\displaystyle \frac{1}{\eta_{1}}\pi+C^{s}_{0}:\D(x,\uu)+
a^{s}_{0}\div_x \uu +a^{s}_{1}p=0.
\end{array} \right\}
\end{equation}
 As before, tensor
$\A^{s}_{0}$, matrices  $C^{s}_{0}, B^{s}_{0}$ and  $B^{s}_{1}$
and constants $a^{s}_{0}$ и $a^{s}_{1}$  are defined by formulas
\eqref{(5.26)}, \eqref{(5.27.1)} - \eqref{(5.27.3)}, in which we
have $\eta_{0}=\eta_{1}$ and $\lambda_{0}=\lambda_{1}$.

Note, that here $\nu_0=0$. Therefore the state equation $ p +\nu_0
p_{*}^{-1}\partial p / \partial t= q$  becomes $p=q$.

 The problem is endowed by the corresponding
homogeneous initial
and boundary conditions.\\

 \qed

\addtocounter{section}{1} \setcounter{equation}{0}
\setcounter{theorem}{0} \setcounter{lemma}{0}
\setcounter{proposition}{0} \setcounter{corollary}{0}
\setcounter{definition}{0} \setcounter{assumption}{0}
\begin{center} \textbf{\S7. Proof of Theorem \ref{theorem2.4}}
\end{center}

\textbf{7.1. Weak and two-scale limits of sequences of
displacement and pressures.}

On the strength of Theorem \ref{theorem2.1}, the sequences
$\{p^\varepsilon\}$, $\{q^\varepsilon\}$, $\{\pi^\varepsilon\}$
and $\{\w^\varepsilon \}$  are uniformly in $\varepsilon$ bounded
in $L^2(\Omega _{T})$. Then there exist a subsequence from
$\{\varepsilon>0\}$ and functions $p$, $\pi$, $q$, and $\w$ such
that as $\varepsilon\searrow 0$
\begin{equation}\label{(7.1.1)}
\w^\varepsilon \rightarrow \w, \quad p^\varepsilon \rightarrow
p,\quad q^\varepsilon \rightarrow q,  \quad \pi^\varepsilon
\rightarrow \pi  \quad \mbox{слабо в } L^2(\Omega_T).
\end{equation}

Moreover, since $\lambda_0,\mu_0>0$ then the bound \eqref{(18)}
imply
\begin{equation} \label{(7.2.2)}
\nabla_x \w^\varepsilon \underset{\varepsilon\searrow
0}{\longrightarrow} \nabla_x \w \quad \mbox{weakly in }
L^2(\Omega_T).
\end{equation}

Due to limiting relations  \eqref{(7.1.1)}, \eqref{(7.2.2)} and
Ngutseng's theorem, there exist one more subsequence from
$\{\varepsilon>0\}$ and 1-periodic in $\y$ functions $P(\x,t,\y)$,
$\Pi(\x,t,\y)$, $Q(\x,t,\y)$, and $\W(\x,t,\y)$ such that the
sequences $\{p^\varepsilon\}$, $\{\pi^\varepsilon\}$,
$\{q^\varepsilon\}$  and  $\{\nabla \w^\varepsilon \}$  two-scale
converge as $\varepsilon \searrow 0$ respectively to $P$, $\Pi$,
$Q $  and  $\nabla _{x}\w +\nabla _y \W $. \\

\textbf{7.2. Micro- and macroscopic equations.}

In the present section we do not consider functions of time $t$,
which re-normalize pressures. As we have shown before, we can
ignore all functions of time.
\begin{lemma} \label{lemma7.1}
Two-scale limits of the sequences $\{p^\varepsilon\}$,
$\{\pi^\varepsilon\}$, $\{q^\varepsilon\}$  and $\{\nabla
\w^\varepsilon\}$  satisfy in $Y_{T}=Y\times (0,T)$ the following
relations
\begin{equation}\label{(7.1)}
    \frac{1}{\eta_{0}}\Pi +(1- \chi)(\div_x \w + \div_y\W) = 0 ;
\end{equation}
\begin{equation}\label{(7.2)}
\frac{1}{p_{*}}P+ \chi (\div_x \w+\div_y\W) =0, \quad Q=P
+\frac{\nu_0}{p_*} \frac{\partial P}{\partial t};
\end{equation}
\begin{eqnarray}\label{(7.3)}
\div_y \bigl(\chi \mu_0 (\D(x,\frac{\partial \w}{\partial
t})+\D(y,\frac{\partial \W}{\partial t})) +
(1-\chi)\lambda_0 (\D(x,\w)+\D(y,\W))\bigr)\\
\nonumber -\nabla _{y}(Q+\Pi )=0.
\end{eqnarray}
\end{lemma}
\begin{lemma} \label{lemma7.2}
The weak limits  $p$, $\pi$, $q$ и  $\w$ satisfy in  $\Omega_T$
the following system of macroscopic equations:
\begin{equation}\label{(7.4)}
   \frac{1}{\eta_{0}}\pi + (1- m)\div_x \w
   +\langle \div_y\W\rangle _{Y_{s}}= 0;
\end{equation}
\begin{equation}\label{(7.5)}
\frac{1}{p_{*}}p+ m\div_x \w+\langle  \div_y\W \rangle _{Y_{f}}
=0,
 \quad
 q=p +\frac{\nu_0}{p_*}\frac{\partial p}{\partial t} ;
\end{equation}
\begin{eqnarray} \label{(7.6)}
&& \tau_0 \hat{\rho}\frac{\partial ^2\w}{\partial t^2} = \div_x
\bigl(\mu_0( m\D(x,\frac{\partial \w}{\partial t})+
 \langle \D(y,\frac{\partial \W}{\partial t})\rangle _{Y_{f}})+\\
 &&\lambda_0 ((1-m)\D(x,\w)+
 \langle \D(y,\W)\rangle _{Y_{s}})
  -(q+\pi)\I\bigr)+\hat{\rho}\F.\nonumber
\end{eqnarray}
\end{lemma}
Proofs of these statements are the same as in lemmas
\ref{lemma5.1} -- \ref{lemma5.3}.

\noindent \textbf{7.3. Homogenized  equations.}
\begin{lemma} \label{lemma7.3}
Weak limits $p$, $\pi$, $q$ and  $\w$ satisfy in  $\Omega_T$ the
following system of homogenized equations:
\begin{eqnarray} \label{(7.7)}
&& \tau_0 \hat{\rho}\frac{\partial ^2\w}{\partial t^2} + \nabla (q+\pi)-\hat{\rho}\F=\\
&&\div_x \bigl(\A_{2}:
\D(x,\frac{\partial \w}{\partial t})+\A_{3}: \D(x,\w)+B_{4}\div_x \w +\nonumber\\
 &&\int_0^t \bigl(\A_{4}(t-\tau ):\D(x,\w(\x,\tau )) +
 B_{5}(t-\tau )\div_x \w(\x,\tau )\bigr)d\tau
  \bigr), \nonumber
\end{eqnarray}
\begin{eqnarray}\label{(7.8)}
&&\frac{1}{p_{*}}p+ m\div_x \w=\\
&&-\int_0^t \bigl(C_{2}(t-\tau ):\D(x,\w(\x,\tau )) +
 a_{2}(t-\tau )\div_x \w(\x,\tau )\bigr)d\tau \nonumber,
\end{eqnarray}

\begin{eqnarray}\label{(7.8.1)}
&& \frac{1}{\eta_{0}}\pi + (1- m)\div_x \w= \\
&&-\int_0^t \bigl(C_{3}(t-\tau ):\D(x,\w(\x,\tau )) +
 a_{3}(t-\tau )\div_x \w(\x,\tau )\bigr)d\tau \nonumber,
\end{eqnarray}
\begin{equation}\label{7.8.2}
 q=p +\frac{\nu_0}{p_*} \frac{\partial p}{\partial t}.
\end{equation}
Here $\A_{2}$, $\A_{3}$ and $\A_{4}$  -- fourth-rank tensors,
$B_{4}$, $B_{5}$,  $C_{2}$  and  $C_{3}$ - matrices and  $a_{2}$
and $a_{3}$- scalars. The exact expressions for these objects are
given below by formulas \eqref{(7.20)}--\eqref{(7.25)}.
\end{lemma}

\begin{proof}

 Let
$$Z(\x,t)=\mu_0\D(x,\frac{\partial \w}{\partial t})-\lambda_0\D(x,\w), \quad Z_{ij}=
\textbf{e}_{i}\cdot (Z\cdot \textbf{e}_{j}), \quad z(\x,t)=\div_x
\w.$$ As usual we look for the solution of the system of
microscopic equations \eqref{(7.1)}--\eqref{(7.3)} in the form
\begin{equation}
\label{(7.9)} \W=\int_0^t\bigl[\W^{0}(y,t-\tau)z(\x,\tau) +
\sum_{i,j=1}^{3}\W^{ij}(\y,t-\tau) Z_{ij}(\x,\tau) \bigr]d\tau,
\end{equation}
\begin{equation}\label{(7.8.3)}
P= \chi \int_0^t\bigl[P^{0}(\y,t-\tau)z(\x,\tau) +
\sum_{i,j=1}^{3}P^{ij}(\y,t-\tau) Z_{ij}(\x,\tau) \bigr]d\tau ,
\end{equation}
\begin{eqnarray}\label{(7.9)}
&&Q=
\chi (Q_{0}(y)\cdot z(\x,t) + \sum_{i,j=1}^{3}Q_{0}^{ij}(\y) \cdot Z_{ij}(\x,t)+\\
&&\int_0^t\bigl[Q^{0}(\y,t-\tau)z(\x,\tau) +
\sum_{i,j=1}^{3}Q^{ij}(\y,t-\tau) Z_{ij}(\x,\tau) \bigr]d\tau)
\nonumber,
\end{eqnarray}
\begin{equation}\label{(7.10)}
\Pi = (1-\chi) \int_0^t\bigl[\Pi^{0}(\y,t-\tau)z(\x,\tau) +
\sum_{i,j=1}^{3}\Pi^{ij}(\y,t-\tau) Z_{ij}(\x,\tau) \bigr]d\tau,
\end{equation}

where 1-periodic in $\y$ functions $\W^{0}$, $\W^{ij}$,  $P^{0}$,
$P^{ij}$, $Q_{0}$,  $Q^{0}$, $Q^{ij}$, $Q_{0}^{ij}$, $\Pi^{0}$,
$\Pi^{ij}$  satisfy the following periodic initial-boundary value
problems
 in the elementary cell $Y$:

\textbf{ Problem ($I$)}
\begin{eqnarray}
\label{(7.11)} && \displaystyle \div_y \bigl(\chi(\mu_0
\D(y,\frac{\partial \W^{ij}}{\partial t}) +\nonumber \\
&&(1-\chi) (\lambda_0 \D(y,\W^{ij}) - ((1-\chi)\Pi^{ij}+\chi Q^{ij}) \I)\bigr)=0;\\
 &&\frac{1}{p_{*}}P^{ij} + \chi \div_y\W^{ij} = 0, \quad
Q^{ij}=P^{ij} +\frac{\nu_0}{p_*} \frac{\partial P^{ij}}{\partial t},\nonumber \\
&&\frac{1}{\eta_{0}}\Pi^{ij}+(1-\chi) \div_y\W^{ij}=0, \quad
 \W^{ij}(\y,0)= \W^{ij}_{0}(\y);\label{(7.13)}\\
&&\displaystyle \div_y \bigl(\chi(\mu_0 \D(y,\W^{ij}_{0})+ J^{ij}
 - Q^{ij}_{0} \I)\bigr)=0,\label{(7.14)}\\
 && \chi(Q^{ij}_{0}+\nu_0 \div_y\W_{0}^{ij})=0. \label{(7.12)}
\end{eqnarray}

\textbf{ Problem ($II$)}
\begin{eqnarray}
\label{(7.15)} && \displaystyle \div_y \bigl(\chi(\mu_0
\D(y,\frac{\partial \W^{0}}{\partial t}) +\nonumber \\
&&(1-\chi) (\lambda_0
\D(y,\W^{0}) - ((1-\chi)\Pi^{0}+\chi Q^{0}) \I)\bigr)=0;\\
 &&\chi(\frac{1}{p_{*}}P^{0} +  \div_y\W^{0}+1) =0, \quad
Q^{0}=P^{0} +\frac{\nu_0}{p_*} \frac{\partial P^{0}}{\partial t}; \label{(7.16)}\\
&&(1-\chi)(\frac{1}{\eta_{0}}\Pi^{0}+ \div_y\W^{0}+1)=0 ;\label{(7.17)}\\
&& \W^{0}(\y,0)= \W^{0}_{0}(\y), \quad \div_y \bigl(\chi (\mu_0
\D(y,\W^{0}_{0})
 - Q_{0} \I)\bigr)=0,\label{(7.18)}\\
&& \chi(Q_{0}+\nu_0 (\div_y\W_{0}^{0}+1))=0.\label{(7.19)}
\end{eqnarray}

Then
\begin{eqnarray}\label{(7.20)}
&&\A_{2}=\mu_0 m \sum_{i,j=1}^{3} J^{ij}\otimes J^{ij} + \mu_0
\A_{0}^{f}, \nonumber \\
&&\A_{0}^{f}=\sum_{i,j=1}^{3}\langle (\mu_0
\D(y,\W^{ij}_{0})\rangle _{Y_{f}} \otimes J^{ij};
\end{eqnarray}
\begin{eqnarray}\label{(7.21)}
&&\A_{3}=\lambda_0(1-m) \sum_{i,j=1}^{3} J^{ij}\otimes J^{ij}-
\lambda_0 \A_{0}^{f} + \mu_0 \A_{1}^{f}(0), \nonumber \\
&&\A_{4}(t)=\mu_0\frac{d}{dt}\A_{1}^{f}(t) - \lambda_0
\A_{1}^{f}(t);
\end{eqnarray}

\begin{equation}\label{(7.22)}
 \A_{1}^{f}(t)=\sum_{i,j=1}^{3}
(\langle \mu_0 \D(y,\frac{\partial \W^{ij}}{\partial
t}(\y,t))\rangle _{Y_{f}}+\langle \lambda_0
\D(y,\W^{ij}(\y,t))\rangle _{Y_{s}})  \otimes J^{ij};
\end{equation}

\begin{equation}\label{(7.23)}
 B_{5}(t)= \mu_0 \langle \D(y,\frac{\partial \W^{0}}{\partial t}\rangle _{Y_{f}}+
  \lambda_0 \langle \D(y,\W^{0}(\y,t))\rangle _{Y_{s}};
\end{equation}

\begin{equation}\label{(7.24)}
 C_{2}(t)=-C_{3}(t)=\sum_{i,j=1}^{3}\langle
 \div_y\W^{ij}(\y,t)\rangle _{Y_{f}} J^{ij};
\end{equation}

\begin{equation}\label{(7.25)}
 a_{2}(t)=-a_{3}(t)=\langle
 \div_y\W^{0}(\y,t)\rangle _{Y_{f}}, \quad
  B_{4}=\mu_0 \langle \D(y,\W^{0}_{0}(\y))\rangle _{Y_{f}}.
\end{equation}

\end{proof}

\begin{lemma} \label{lemma7.4}
Tensors $\A_{2}$-- $\A_{4}$, matrices $B_{4}$, $B_{5}$,  $C_{2}$
and $C_{3}$ and scalars  $a_{2}$ and  $a_{3}$ are well-defined and
infinitely smooth in time.

If a porous space is connected, then the symmetric tensor $\A_{2}$
is strictly positively defined. For the case of disconnected
porous space (isolated pores) $\A_{2}=0$ and the tensor $\A_{2}$
becomes strictly positively defined.
\end{lemma}

All these objects are well-defined if \textbf{ Problem ($I$)} and
\textbf{ Problem($II$)} are well-posed. The solvability of above
mentioned problems and smoothness with respect to time follow, due
to linearity, from the standard a'priory estimates.

The symmetry of  $\A_{2}$  proves in the same way as the symmetry
of $\A_{0}$. If the porous space is disconnected, then the problem
\eqref{(7.14)} has a unique solution linear in $\y$, such that

\begin{equation}\label{(7.26)}
    \chi(\D(y,\W^{ij}_{0})+ J^{ij})=0.
\end{equation}
The last equality implies  $\A_{2}=0.$

In this case the tensor $\A_{3}$ becomes strictly positively
defined. Indeed

$$\A_{3}=\lambda_0 \sum_{i,j=1}^{3} J^{ij}\otimes J^{ij}
 + \mu_0 \A_{1}^{f}(0)=$$
 $$\lambda_0 \sum_{i,j=1}^{3} J^{ij}\otimes J^{ij}+
 \sum_{i,j=1}^{3}
\langle \chi\mu_0 \D(y,\frac{\partial \W^{ij}}{\partial t}(\y,0))+
\frac{\lambda_0}{\mu_0} J^{ij}\rangle _{Y} \otimes J^{ij}.$$ On
the other hand, coming back to \eqref{(7.11)} at initial time
moment we see that
  $$\langle \chi\mu_0 \D(y,\frac{\partial \W^{ij}}{\partial t}(\y,0)):\D(y,\W^{kl}_{0})\rangle
  _{Y}=$$
  $$-\lambda_0\langle \chi \D(y,\W^{ij}_{0}):\D(y,\W^{kl}_{0})\rangle _{Y}-
  \langle \frac{1}{\eta_{0}}\Pi^{ij}\cdot\Pi^{kl} \rangle _{Y}|_{t=0}.$$
  Moreover, due to \eqref{(7.26)}
$$\langle \chi\mu_0 \D(y,\frac{\partial \W^{ij}}{\partial t}(\y,0)):\D(y,\W^{kl}_{0})\rangle=
-\langle \chi \D(y,\frac{\partial \W^{ij}}{\partial
t}(\y,0)):J^{kl}\rangle _{Y},$$ which proves our statement.

Anvarbek  Meirmanov

Ugra State University, Khanti-Mansiisk, Russia;

Center for Advanced Mathematics and Physics, Electrical and
Mechanical Engineering College, National University of Science and
Technology, Peshawar Road, Rawalpindi, Pakistan;

email: anvarbek@list.ru

\end{document}